\documentclass[12pt,reqno]{amsart}

\usepackage{mathtools,amssymb,amsthm}
\usepackage{bm}
\usepackage{empheq, esint}

\usepackage[a4paper,margin=2.5cm]{geometry}
\usepackage{xcolor}
\usepackage{graphicx}
\usepackage{subcaption}

\usepackage[english]{babel}
\usepackage{csquotes}
\usepackage[backref=page]{hyperref}
\usepackage{orcidlink}

\newtheorem{theorem}{Theorem}[section]
\newtheorem{lemma}[theorem]{Lemma}

\newtheorem{proposition}[theorem]{Proposition}

\theoremstyle{definition}
\newtheorem{definition}[theorem]{Definition}

\theoremstyle{remark}

\numberwithin{equation}{section}

\newcommand{\RR}{\ensuremath{\mathbb{R}}}

\DeclareMathOperator{\supp}{supp}
\DeclareMathOperator{\loc}{loc}

\newcommand{\md}{\mathrm{\,d}}

\newcommand{\xqedhere}[2]{\rlap{\hbox to#1{\hfil\llap{\ensuremath{#2}}}}}
\newcommand{\orcid}[1]{\textsuperscript{\orcidlink{#1}}}

\title[Nonexistence of positive solutions stable outside a compact set]{Nonexistence of positive solutions stable outside a compact set for the Lane--Emden system in a half-space}

\author[Zetao Cheng]{Zetao Cheng}
\address{Department of Mathematics, Hanyang University, Seoul 04763, Republic of Korea}
\email{\href{mailto:chengzt20@hanyang.ac.kr}{\texttt{chengzt20@hanyang.ac.kr}}}

\author[Long-Han Huang]{Long-Han Huang \orcid{0009-0003-1864-2829}}
\address{Department of Mathematical Sciences, Tsinghua University, Beijing 100084, China}
\email{\href{mailto:huanglh22@mails.tsinghua.edu.cn}{\texttt{huanglh22@mails.tsinghua.edu.cn}}}

\author[Wenming Zou]{Wenming Zou \orcid{0000-0002-4178-9568}}
\address{Department of Mathematical Sciences, Tsinghua University, Beijing 100084, China}
\email{\href{mailto:zou-wm@mail.tsinghua.edu.cn}{\texttt{zou-wm@mail.tsinghua.edu.cn}}}

\date{November 6, 2026}

\begin{document}

\begin{abstract}
We prove that, in every dimension and for all positive exponents, the Lane--Emden system in a half-space with homogeneous Dirichlet boundary conditions admits no positive classical solution stable outside a compact set.
\end{abstract}

\maketitle
\tableofcontents

\section{Introduction}

For $n\ge1$, let $\mathbb{R}_+^n=\{x=(x',x_n)\in\RR^{n-1}\times\RR:x_n>0\}$. We consider the Lane--Emden system with homogeneous Dirichlet boundary conditions,
\begin{align}\label{eq:LE}
    \begin{aligned}
        -\Delta u&=v^p,& -\Delta v&=u^q &&\text{in }\mathbb{R}_+^n,\\
        u&=0,& v&=0 &&\text{on }\partial\mathbb{R}_+^n,
    \end{aligned}
\end{align}
where $p,q>0$. By a classical solution, we mean a pair $(u,v)\in[C^2(\mathbb{R}_+^n)\cap C(\overline{\mathbb{R}_+^n})]^2$ satisfying the equations and boundary conditions pointwise. By the strong minimum principle, every nonnegative, nontrivial solution is strictly positive in $\mathbb{R}_+^n$.

Motivated by the notions of stability for cooperative elliptic systems in~\cite{Montenegro2005}, \cite[Definition 1]{Cowan2013}, \cite[Definition 1.1]{ChenDupaigneGhergu2014}, and~\cite[Definition 1.1]{HuangZou2026}, we make the following definition; see also~\cite[Introduction]{MtiriYe2019,DupaigneGherguHajlaoui2025} for stability outside a compact set.

\begin{definition}\label{def:stability}
Let $n\ge1$, $p,q>0$, and let $(u,v)$ be a positive classical solution of~\eqref{eq:LE}. The solution is \emph{stable in an open set $\Omega\subseteq \mathbb{R}_+^n$} if there exist positive functions $\Phi,\Psi\in C^2(\Omega)$ satisfying
\begin{align}\label{eq:Jacobi}
    \begin{aligned}
        -\Delta\Phi&\ge p v^{p-1}\Psi,\\
        -\Delta\Psi&\ge q u^{q-1}\Phi
    \end{aligned}
    \qquad\text{in }\Omega.
\end{align}
The solution is \emph{stable outside a compact set} if there is a compact $K\subseteq\overline{\mathbb{R}_+^n}$ such that it is stable in every open set $\Omega\Subset\mathbb{R}_+^n\setminus K$.
\end{definition}

For the compact set $K$ in Definition~\ref{def:stability}, Lemma~\ref{lem:local-stability} shows that stability outside $K$ is equivalent to the existence of positive functions $\Phi,\Psi\in C^2(\mathbb{R}_+^n\setminus K)$ satisfying~\eqref{eq:Jacobi} throughout $\mathbb{R}_+^n\setminus K$. More generally, by Lemmas~\ref{lem:stability} and~\ref{lem:local-stability}, a positive classical solution $(u,v)$ is stable in an open set $\Omega\subseteq\mathbb{R}_+^n$ if and only if, for every bounded smooth connected open set $D\Subset\Omega$,
\begin{equation*}
\|\sqrt{pv^{p-1}}\,G_D\sqrt{qu^{q-1}}\|_{L^2(D)\to L^2(D)}\le1\,.
\end{equation*}
Here $G_D$ denotes the Dirichlet Green operator on $D$.

Stability in the whole half-space implies stability outside a compact set. Monotonicity in the normal direction also implies stability in the whole half-space: if $u_{x_n},v_{x_n}\ge0$, then differentiating~\eqref{eq:LE} gives the linearized system for $(u_{x_n},v_{x_n})$, and the boundary conditions and the strong minimum principle give $u_{x_n},v_{x_n}>0$. For $n\ge2$ and $p,q>1$, Li and Souplet~\cite[Proposition 1.2]{LiSouplet2025} proved that every positive classical solution bounded on each strip $\mathbb{R}^{n-1}\times(0,R)$, $R>0$, satisfies $u_{x_n},v_{x_n}>0$. Such solutions are therefore stable in the whole half-space.

For $n\ge2$, $p,q>0$, and $pq>1$, Lemma~\ref{lem:global-reduction} shows that system~\eqref{eq:LE} admits a positive classical solution stable outside a compact set if and only if it admits a positive classical solution stable in the whole half-space. The globally stable solution constructed in that lemma may differ from the original solution.

We now state our main nonexistence theorem.

\begin{theorem}\label{thm:main}
For every $n\ge1$ and $p,q>0$, system~\eqref{eq:LE} admits no positive classical solution stable outside a compact set.
\end{theorem}

Theorem~\ref{thm:main} assumes neither boundedness on finite strips nor any growth or energy bounds or monotonicity in the normal direction. It also covers the case $0<\min\{p,q\}<1$.

We next review related Liouville theorems for scalar equations and elliptic systems, focusing on the assumptions concerning boundedness, monotonicity, and stability.

\subsection{Scalar equations}

Half-space Liouville theorems play an important role in blow-up arguments for establishing a priori bounds near the boundary; see~\cite{GidasSpruck1981CPDE,Giga1986,PolacikQuittnerSouplet2007}. For the scalar equation $-\Delta w=w^r$ in $\mathbb{R}_+^n$ with zero Dirichlet boundary data, Gidas and Spruck~\cite[Theorem 1.3 and Remark 2 following its proof]{GidasSpruck1981CPDE} proved nonexistence of positive solutions when $n\ge3$ and $1<r\le(n+2)/(n-2)$, without assuming boundedness. Their result includes the critical exponent. For bounded positive solutions, Dancer~\cite[Theorem 2 and Remarks 2]{Dancer1992} used monotonicity in the normal direction and a limit in one fewer dimension to obtain nonexistence for $1<r<(n+1)/(n-3)$ when $n>3$, and for every $r>1$ when $2\le n\le3$. Using stability, Farina obtained nonexistence for every $r>1$ for bounded positive solutions when $2\le n\le11$~\cite[Theorem 12]{Farina2007}, and for positive solutions stable outside a compact set when $2\le n\le10$, without assuming boundedness~\cite[Theorem 9(b)]{Farina2007}.

Chen, Lin, and Zou~\cite{ChenLinZou2014} proved nonexistence of bounded positive solutions for every $r>1$ by establishing convexity in the normal direction. Dupaigne, Sirakov, and Souplet~\cite[Theorem 1.1 and Corollary 1.2]{DupaigneSirakovSouplet2022} extended this result to solutions bounded on each finite strip or monotone in the normal direction, without further growth assumptions. Dupaigne, Farina, and Petitt~\cite[Theorem 1]{DupaigneFarinaPetitt2023} proved nonexistence under stability outside a compact set in every dimension and for every $r>1$. The equation considered there, $-\Delta w=|w|^{r-1}w$, also allows sign-changing solutions. For $n\ge3$ and $0<r\le1$, Luo~\cite[Theorem 1.1]{Luo2026} proved nonexistence of nonnegative nontrivial solutions without a stability assumption.

\subsection{Elliptic systems}

For $n\ge3$ and $p,q>1$, Birindelli and Mitidieri~\cite[Theorem 5.1 and the argument leading to (5.2)]{BirindelliMitidieri1998} used normal monotonicity to obtain solutions in one fewer dimension from bounded positive Dirichlet solutions of fully coupled cooperative power systems in a half-space. For the Lane--Emden system, de Figueiredo and Sirakov~\cite[Theorems 1.2--1.3 and Sections 3.3--3.4]{deFigueiredoSirakov2005} established strict normal monotonicity and showed that bounded positive solutions converge, as the normal variable tends to infinity, to bounded positive entire solutions in $\RR^{n-1}$. Their limiting argument relies on the boundedness assumption.

Reichel and Zou~\cite[Theorem 5 and its proof in Section 5]{ReichelZou2000} proved one-dimensionality of nonnegative Dirichlet solutions for autonomous cooperative systems with locally Lipschitz nonlinearities satisfying their componentwise subcriticality condition, without assuming boundedness. In particular, for $n\ge3$ and $1<p,q\le(n+2)/(n-2)$, their result excludes positive solutions of~\eqref{eq:LE}, since the one-dimensional equations admit no positive solution on the half-line.

A different connection with entire solutions is given by~\cite[Theorem 4.2]{PolacikQuittnerSouplet2007}. For fixed $n\ge1$ and $p,q>1$, nonexistence of bounded nontrivial nonnegative solutions in $\RR^n$ implies nonexistence of nontrivial nonnegative classical Dirichlet solutions in the half-space of the same dimension.

Stability methods also yield nonexistence for bounded solutions of the half-space system. Cowan~\cite[Theorem 3]{Cowan2013} obtained this conclusion for $2\le n\le11$ and $\min(p,q)\ge2$. Chen, Lin, and Zou~\cite[Corollary 1.1]{ChenLinZou2014} proved that every bounded nonnegative solution is zero for all $p,q>1$, in every dimension.

\subsection{Comparison with recent work}

\emph{While completing this paper, we became aware of the recent preprint by Li and Souplet~\cite{LiSouplet2026} on monotone solutions of the Lane--Emden system in a half-space. The present work was developed independently of that preprint.}

In their earlier paper, Li and Souplet~\cite[Theorem 1.1]{LiSouplet2025} proved nonexistence for the half-space system with $p,q>1$ under boundedness on each finite strip. Their proof combines moving planes, logarithmic gradient estimates for the normal derivatives, and convexity in the normal direction. Their recent preprint~\cite[Theorem 1.1]{LiSouplet2026} establishes nonexistence for $n\ge2$ and $p,q\ge1$ under the assumption $u_{x_n},v_{x_n}\ge0$.

As noted above, the monotone solutions considered in~\cite[Theorem 1.1]{LiSouplet2026} are globally stable; see also~\cite[Proposition 2.5]{LiSouplet2026}. Theorem~\ref{thm:main} requires stability only outside a compact set and also allows positive exponents below one.

For $pq>1$, Li and Souplet's proof uses a comparison on finite strips in place of normal convexity. They form affine interpolants from the normal limits and project the positive parts of the differences between the solution and these interpolants onto the first sine mode. A strict contraction estimate for a Bessel resolvent, obtained from the positive normal derivatives through a two-weight Schur test, then forces the projections to vanish. The case $p=q=1$ follows from a separate a priori estimate that does not require normal monotonicity~\cite[Remark 2.3]{LiSouplet2026}.

Our strip comparison uses normal monotonicity only for auxiliary minimal solutions. For each height $S>0$, we construct the minimal pair $(M_S,N_S)$ with the original top boundary data and use its monotonicity on the lower half of the strip. With $\rho=\min(q,1)$, we average $M_S^\rho$ and $N_S$ against a polynomial weight whose value and first derivative vanish at both endpoints. The $1/\rho$ power of the first average and the second average form a nonlinear supersolution in one fewer variable. Combining this supersolution with an explicit quadratic barrier gives a supersolution on a lower strip that dominates the Picard iterates defining $(M_S,N_S)$. This comparison uses neither the normal limit of the original solution nor a projection of an affine interpolation error.

When $\rho<1$, the positive gradient term in the equation for $M_S^\rho$ compensates for the negative gradient term introduced by taking the $1/\rho$ power of its average. H\"older's inequality then gives the nonlinear supersolution inequalities, allowing the averaging argument to include a sublinear exponent.

The two arguments share the use of finite strips, normal averaging, stability inequalities, and Green operators. We obtain the Green contraction from a fourth-order Picone inequality on interior domains, without a two-weight Schur test. The factor $(pq)^{-1/2}$ also follows from~\eqref{eq:green-contraction} by removing $p$ and $q$ from its two weights, so the numerical constants alone do not distinguish the Green bounds.

For the reaction-energy estimate, we follow~\cite[proof of Lemma 4]{Cowan2013} and combine a power test in one component with the pointwise comparison between components. The two-component tests in~\cite[Lemma 2.6]{LiSouplet2026} use the same stability inequality. In our estimate, the component comparison makes the coefficient of the reaction integral on the left strictly larger than the coefficient on the right, so a single component test suffices to bound that integral by cutoff terms. A separate weighted estimate treats the case of a sublinear exponent.

Exponential weights also occur in both proofs. In~\cite[Lemma 5.6 and the proof of Proposition 5.7]{LiSouplet2026}, they preserve strict contraction of the weighted Bessel resolvent and force the projected interpolation errors to vanish. Here exponential conjugation gives tangential localization estimates for finite powers of the full linearized operator $Q_h$ defined in~\eqref{eq:Q-definition}, whose coefficients come from the original solution. We use only the bound $\|Q_h\|\le1$ and obtain a contradiction from nonlinear mass estimates at finitely many fixed strip heights.

\subsection{Proof outline}

The cases $n=1$ and $pq\le1$ follow from concavity and Green comparison, respectively. We therefore consider $n\ge2$ and $pq>1$, and interchange the components if necessary so that $p\ge q$. The proof has three main parts.

First, Lemma~\ref{lem:global-reduction} constructs a globally stable positive solution from a positive solution stable outside a compact set. We construct minimal solutions on expanding half-balls and show that a subsequential limit is positive. Convergence to zero would imply $u^{q+1},v^{p+1}\in L^1(\mathbb{R}_+^n)$ for the original solution. The normal-translation identity would then force the integral of the strictly positive boundary flux $u_{x_n} v_{x_n}$ to vanish, contradicting Hopf's lemma. When $0<q<1$, the inverse-power difference in~\eqref{eq:inverse-power-difference} preserves the coefficient $qu^{q-1}$ of the original linearization and produces a nonnegative remainder. This allows the reduction to include sublinear exponents and, for $n\ge3$, exponents on the Sobolev critical hyperbola $n/(p+1)+n/(q+1)=n-2$.

Second, Lemmas~\ref{lem:weighted-trace} and~\ref{lem:strip-decay} give decay estimates for the minimal strip pair without assuming stability of the original solution. Averaging powers of $(M_S,N_S)$ and using a quadratic barrier gives a supersolution that dominates the Picard iterates defining this pair. A local doubling argument and a boundary gradient estimate then yield constants $c,C>0$, depending only on $n,p,q$, such that for every $S>0$, $x'\in\RR^{n-1}$, and $0<x_n\le cS$,
\begin{equation*}
M_S(x',x_n)\le C x_n S^{-\alpha-1},\qquad N_S(x',x_n)\le C x_n S^{-\beta-1}\,,
\end{equation*}
where $\alpha,\beta$ are defined in~\eqref{eq:scaling-exponents}. These bounds are uniform in $x'$ and independent of the original solution. They concern the auxiliary pair $(M_S,N_S)$ and imply neither decay nor monotonicity of $(u,v)$. Their uniformity allows the sources in the final argument to move arbitrarily far in the tangential directions.

Third, we combine strip decay with estimates for the positive self-adjoint contraction $Q_h$ on $L^2(\mu)$, where $\md\mu=pv^{p+1}\md x$. We use the bound $\|Q_h\|\le1$ without requiring a spectral gap below one. Strict convexity of the composition of the nonlinear Green maps gives quadratic lower bounds for the masses of finite sums of iterates. To repeat this estimate at increasing strip heights, Lemma~\ref{lem:source-energy} bounds the energy of each new source by finite energy differences associated with a single compactly supported source. For sums of between $N$ and $2N$ iterates, each increase to a suitable fixed strip height multiplies the mass lower bound by $cN^2$, while the $L^2(\mu)$ norm remains $O(\sqrt N)$ at every fixed height. Lemma~\ref{lem:localization} shows that the $\mu$-mass outside a horizontal radius of order $\sqrt{N\log N}$ tends to zero. After finitely many height increases, the mass lower bound exceeds the upper bound obtained from Cauchy--Schwarz and the polynomial growth of $\mu$. All strip heights and constants are chosen independently of $N$.

Section~\ref{sec:preliminaries} proves the source and stability estimates. Section~\ref{sec:global} constructs a globally stable solution from an exterior-stable one. Section~\ref{sec:strips} proves the uniform decay estimate for minimal strip solutions. Section~\ref{sec:operators} establishes the energy and operator bounds and completes the proof.

\section{Source estimates and the meaning of stability}\label{sec:preliminaries}

We first derive estimates for positive solutions of~\eqref{eq:LE} without assuming stability. These estimates will be used in Green comparisons on unbounded strips and in integrations by parts near the boundary.

For $R>0$ and $\xi\in\RR^{n-1}$, we write $B_R^+=\{x\in\RR^{n}:|x|<R\}\cap\mathbb{R}_+^n$ and $B'_R(\xi)=\{x'\in\RR^{n-1}:|x'-\xi|<R\}$.

For an open set $D\subseteq\mathbb{R}_+^n$, let $G_D$ denote the minimal Dirichlet Green operator associated with $-\Delta$; see~\cite[Section 4.1]{ArmitageGardiner2001} and~\cite[Section 5.7]{Grigoryan2006}. Choose smooth bounded open sets $D_j$ such that $D_j\Subset D_{j+1}\Subset D$ and $\bigcup_jD_j=D$. Let $G_{D_j}(x,y)$ be the Dirichlet Green kernel on $D_j$, extended by zero to $D\times D$. By domain monotonicity, for distinct $x,y\in D$,
\begin{equation*}
G_{D_j}(x,y)\nearrow G_D(x,y)\qquad\text{as }j\to\infty\,.
\end{equation*}
The limit is independent of the chosen exhaustion. It is nonnegative and bounded above by the half-space Green kernel, and hence finite for $x\ne y$.

For every nonnegative measurable function $f$ on $D$ and every $x\in D$, define
\begin{align*}
G_Df(x)&:=\int_D G_D(x,y)f(y)\md y\\
&=\lim_{j\to\infty}\int_{D_j}G_{D_j}(x,y)f(y)\md y\,.
\end{align*}
The equality follows from monotone convergence. These are nonnegative Lebesgue integrals and may take the value $+\infty$; finiteness will be established whenever it is needed.

When $pq>1$, set
\begin{equation}\label{eq:scaling-exponents}
\alpha=\frac{2(p+1)}{pq-1},\qquad \beta=\frac{2(q+1)}{pq-1}.
\end{equation}
Thus $\alpha+2=p\beta$, $\beta+2=q\alpha$, and $\alpha(q+1)=\beta(p+1)=\alpha+\beta+2$. These identities give the scaling of~\eqref{eq:LE}.

Let $(u,v)$ be a positive classical solution of~\eqref{eq:LE}. Local boundary regularity will be used repeatedly. On every compact flat boundary patch, the sources are bounded, so boundary Poisson estimates give $u,v\in C^{1,\gamma}$ for $\gamma<1$. The powers are then locally H\"older continuous, and $u,v\in C^{2,\gamma}$ for any sufficiently small $\gamma>0$. Hopf's lemma gives $c x_n\le u,v\le C x_n$ on a sufficiently small compact boundary neighborhood; these constants may depend on the solution and the patch.

For every compact $K'\subseteq\partial\mathbb{R}_+^n$, there exist constants $c,C,\varepsilon>0$, depending on the solution and $K'$, such that
\begin{equation}\label{eq:local-boundary-orders}
c x_n\le u(x',x_n),v(x',x_n)\le C x_n
\qquad((x',0)\in K',\ 0<x_n<\varepsilon)\,.
\end{equation}
Together with the preceding boundary regularity, these local bounds will justify integrations by parts at the flat boundary.

We next establish integral estimates that are uniform in the height and tangential position. In Lemma~\ref{lem:source}, the horizontal-slice estimates yield volume bounds for $u,v$ and weighted integral bounds for the sources $v^p,u^q$. The volume bounds also exclude nonnegative harmonic remainders in the Green representations of $u$ and $v$, giving the identities in part~\textup{(iii)}.

The slice estimates will be used in the comparison and averaging arguments on unbounded strips. The weighted source estimates will control the terms at infinity in Section~\ref{sec:global}, where the Green representation will also give $\int_{\mathbb R_+^n}u^{q+1}=\int_{\mathbb R_+^n}v^{p+1}$. Unless otherwise stated, constants in Lemma~\ref{lem:source} depend only on $n,p,q$.

\begin{lemma}\label{lem:source}
Let $n\ge2$, $p,q>0$, $pq>1$, and let $(u,v)$ be a positive classical solution of~\eqref{eq:LE}.
\begin{enumerate}
\item[\textup{(i)}] For every $h,R>0$ and $\xi\in\RR^{n-1}$,
\begin{align}\label{eq:horizontal-L1}
    \begin{aligned}
        \int_{B'_R(\xi)}u(x',h)\md x'&\le C R^{n-1-\alpha},\\
        \int_{B'_R(\xi)}v(x',h)\md x'&\le C R^{n-1-\beta}.
    \end{aligned}
\end{align}
\item[\textup{(ii)}] In addition, for all $R>0$,
\begin{align}\label{eq:source-moments}
    \begin{aligned}
        \int_{B_R^+}u&\le CR^{n-\alpha},&
        \int_{B_R^+}x_n v^p&\le CR^{n-1-\alpha},\\
        \int_{B_R^+}v&\le CR^{n-\beta},&
        \int_{B_R^+}x_n u^q&\le CR^{n-1-\beta}.
    \end{aligned}
\end{align}
\item[\textup{(iii)}] Finally, $u=G_{\mathbb{R}_+^n}(v^p)$ and $v=G_{\mathbb{R}_+^n}(u^q)$.
\end{enumerate}
\end{lemma}

\begin{proof}
We first use the same-scale Green comparison for systems; compare~\cite[Sections 6.1--6.2]{ArmstrongSirakov2011}. If $B_{8r}(z)\subseteq \mathbb{R}_+^n$, the lower bound for the Green kernel of $B_{2r}(z)$ on $B_r(z)\times B_r(z)$ gives
\[
u(x)\ge c r^{2-n}\int_{B_r(z)}v^p,\qquad
v(x)\ge c r^{2-n}\int_{B_r(z)}u^q\qquad(x\in B_r(z)).
\]
For $n=2$ the factor $r^{2-n}$ is one. Taking positive powers, integrating on the same ball, and eliminating the two integrals yields
\begin{align}\label{eq:interior-source}
    \begin{aligned}
        \int_{B_r(z)}v^p&\le Cr^{n-\alpha-2},\\
        \int_{B_r(z)}u^q&\le Cr^{n-\beta-2}.
    \end{aligned}
\end{align}
This argument uses $pq>1$ but does not require either power to be convex.

To estimate $u$ itself, decompose it in $B_{2r}(z)$ as the sum of its Green potential and a nonnegative harmonic function,
\begin{equation}\label{eq:local-green-decomposition}
 u=G_{B_{2r}(z)}(v^p)+h_0,\qquad h_0\ge0,\quad\Delta h_0=0.
\end{equation} The integral of the potential over $B_r(z)$ is at most $Cr^2\int_{B_{2r}(z)}v^p$, by the torsion function of the ball. For the harmonic part, Harnack's inequality and~\eqref{eq:interior-source} give
\[
\sup_{B_r(z)}h_0\le C\inf_{B_r(z)}h_0
\le C\inf_{B_r(z)}u\le C\left(\fint_{B_r(z)}u^q\right)^{1/q}
\le Cr^{-\alpha}.
\]
Adjusting the fixed radius ratios if necessary, we obtain $\int_{B_r(z)}u\le Cr^{n-\alpha}$. The other component satisfies the analogous estimate.

Fix the plane $x_n=h$. Positivity and superharmonicity imply
\[
u(x',x_n)\ge\int_{B'_R(\xi)}P_{x_n-h}(x'-y')u(y',h)\md y',\qquad x_n>h,
\]
where $P_s$ is the Poisson kernel of the upper half-space. To justify this comparison, first truncate the boundary data to a bounded compactly supported function, compare in bounded half-balls, and then increase the domain and the data. No comparison at infinity is assumed. On $B_{R/16}((\xi,h+R))$ the kernel in this integral is bounded below by $cR^{1-n}$. Integrating there and applying the interior estimate just proved gives~\eqref{eq:horizontal-L1}. Its constants are independent of the height. Integration in $0<x_n<R$ gives the unweighted estimates in~\eqref{eq:source-moments}. Testing the first equation against $x_n\chi(x/R)$, with a fixed $0\le\chi\in C_c^\infty(B_2)$ equal to one on $B_1$, gives
\[
\int_{B_R^+}x_n v^p\le \frac C R\int_{B_{2R}^+}u\le CR^{n-1-\alpha}.
\]
Both functions in the boundary term vanish on the flat boundary; the local regularity preceding Lemma~\ref{lem:source} justifies the integration by parts. The second equation gives the remaining estimate.

Positive Green comparison gives $0\le G_{\mathbb{R}_+^n}(v^p)\le u$. The increasing Green potentials solve the Poisson equation in distributions, so $w=u-G_{\mathbb{R}_+^n}(v^p)$ is harmonic in $\mathbb{R}_+^n$. Since $0\le w\le u$, it extends continuously by zero to the boundary. Its zero extension is subharmonic on $\RR^n$: for each $\varepsilon>0$, the positive part $(w-\varepsilon)_+$ vanishes near each compact boundary patch and can be extended by zero; then let $\varepsilon\downarrow0$ locally in $L^1$. For each fixed $x\in \mathbb{R}_+^n$ and sufficiently large $R$, the mean-value inequality and~\eqref{eq:source-moments} give
\[
w(x)\le CR^{-n}\int_{B_R(x)\cap \mathbb{R}_+^n}w
\le CR^{-n}\int_{B_{2R}^+}u
\le CR^{-\alpha}\longrightarrow0.
\]
Thus $u=G_{\mathbb{R}_+^n}(v^p)$. The same argument for $v-G_{\mathbb{R}_+^n}(u^q)$ uses $\beta>0$ and gives the second representation. This also proves the claim in dimension two.
\end{proof}

\begin{lemma}\label{lem:comparison}
If $p\ge q>0$ and $pq>1$, every positive solution in Lemma~\ref{lem:source} satisfies
\begin{align}\label{eq:component-comparison}
    \frac{v^{p+1}}{p+1}\le\frac{u^{q+1}}{q+1}.
\end{align}
The same inequality holds for a positive solution in a bounded half-ball whose boundary data satisfy this inequality.
\end{lemma}

\begin{proof}
We use the component comparison calculation of~\cite[Lemma 2.7]{Souplet2009}; the unbounded half-space step below also covers $q<1$. Let $\sigma=(q+1)/(p+1)\le1$ and $\ell=\sigma^{-1/(p+1)}$. On the set where $v>\ell u^\sigma$,
\[
-\Delta(v-\ell u^\sigma)
\le u^q-\ell\sigma u^{\sigma-1}v^p\le0.
\]
Consequently $w=(v-\ell u^\sigma)_+$, extended by zero across the flat boundary, is subharmonic on $\RR^n$. For $\sigma<1$ it vanishes near each compact boundary patch by~\eqref{eq:local-boundary-orders}; for $\sigma=1$ the zero-trace positive-part extension gives the same distributional assertion. Since $0\le w\le v$, the mean-value inequality and~\eqref{eq:source-moments} imply $w(x)\le CR^{-\beta}$ for any fixed $x$ and all sufficiently large $R$. Hence $w=0$. In a bounded half-ball the same calculation uses the ordinary maximum principle and the prescribed boundary inequality. The related half-space comparison in~\cite[Proposition 4.1]{LiSouplet2025} is stated under that paper's standing assumption $p,q>1$.
\end{proof}

A positive classical supersolution in $\mathbb R_+^n$ is a pair $(u,v)\in[C^2(\mathbb R_+^n)]^2$ with $u,v>0$ satisfying
\begin{equation}\label{eq:positive-supersolution}
-\Delta u\ge v^p,\qquad -\Delta v\ge u^q\quad\text{in }\mathbb R_+^n.
\end{equation}

\begin{lemma}\label{lem:small-product}
Let $n\ge2$ and $p,q>0$ with $pq\le1$. Then no positive classical supersolution exists in $\mathbb R_+^n$, even without any boundary condition.
\end{lemma}

\begin{proof}
We use the same-scale comparison and algebraic elimination method for systems in~\cite[Sections 6.1--6.2, especially (6.2) and (6.6)]{ArmstrongSirakov2011}, here applied to interior balls receding into the half-space. In $B_R(3Re_n)$ set $A_R=\int_{B_R(3Re_n)}v^p$ and $B_R=\int_{B_R(3Re_n)}u^q$. Green comparison in the concentric ball of radius $2R$ gives
\[
A_R\ge cR^{n-(n-2)p}B_R^p,\qquad
B_R\ge cR^{n-(n-2)q}A_R^q.
\]
Thus $A_R^{1-pq}\ge cR^{n(1-pq)+2p(q+1)}$. If $pq=1$, this is impossible as $R\to\infty$. If $pq<1$, it gives $A_R\ge cR^{n+2p(q+1)/(1-pq)}$. On the other hand, comparison with the half-space Green potential of $v^p\mathbf1_{B_R(3Re_n)}$ at $e_n$ gives $u(e_n)\ge cR^{1-n}A_R$ for $R\ge2$. These bounds are incompatible. The last lower bound also holds for the logarithmic half-space kernel when $n=2$.
\end{proof}

The fourth-order inequality below follows from the full cooperative stability condition. It is essential when one exponent is sublinear. The geometric-mean inequality alone will suffice for the other energy estimate.

\begin{lemma}\label{lem:stability}
Let $n\ge1$ and $p,q>0$, and let $(u,v)$ be a positive classical solution of~\eqref{eq:LE} that is stable in an open set $D\subseteq\mathbb{R}_+^n$. Then the following statements hold.
\begin{enumerate}
\item[\textup{(i)}] The weighted Green operator satisfies
\begin{equation}\label{eq:green-contraction}
\bigl\|\sqrt{pv^{p-1}}\,G_D\sqrt{qu^{q-1}}\bigr\|_{L^2(D)\to L^2(D)}\le1\,.
\end{equation}

\item[\textup{(ii)}] For every $\zeta\in H_0^1(D)$ with compact support in $\overline D$,
\begin{equation}\label{eq:Cowan}
\sqrt{pq}\int_D u^{(q-1)/2}v^{(p-1)/2}\zeta^2
\le\int_D|\nabla\zeta|^2\,.
\end{equation}

\item[\textup{(iii)}] If $D$ is bounded, then every $\phi\in H^2(D)\cap H_0^1(D)$ satisfying $u^{(1-q)/2}\Delta\phi\in L^2(D)$ obeys
\begin{equation}\label{eq:Navier}
pq\int_D v^{p-1}\phi^2
\le\int_D u^{1-q}(\Delta\phi)^2\,.
\end{equation}
\end{enumerate}
\end{lemma}

\begin{proof}
We use the standard Picone argument. Let $\Omega\Subset D$ be smooth. Restrict the Jacobi supersolutions satisfying~\eqref{eq:Jacobi} on $D$ to a neighborhood of $\overline\Omega$. For $\phi\in C^2(\overline\Omega)$ with $\phi=0$ on $\partial\Omega$, both $\phi^2/\Psi$ and its normal derivative vanish on $\partial\Omega$. With $E=\{\Delta(\phi^2/\Psi)<0\}$, two integrations by parts and the two supersolution inequalities give
\begin{align*}
p\int_\Omega v^{p-1}\phi^2
&\le\int_\Omega(-\Delta\Phi)\frac{\phi^2}{\Psi}=\int_\Omega\Phi\left[-\Delta\left(\frac{\phi^2}{\Psi}\right)\right]\\
&\le\int_E\Phi\left[-\Delta\left(\frac{\phi^2}{\Psi}\right)\right]\le\frac1q\int_E u^{1-q}\Delta\Psi\,
\Delta\left(\frac{\phi^2}{\Psi}\right)\,.
\end{align*}
The restriction to $E$ is necessary for the last inequality. Since $\Delta\Psi\le0$, the pointwise identity
\begin{align*}
(\Delta\phi)^2-\Delta\Psi\,\Delta\left(\frac{\phi^2}{\Psi}\right)
&=\left(\Delta\phi-\frac{\phi}{\Psi}\Delta\Psi\right)^2
-\frac{2\Delta\Psi}{\Psi}
\left|\nabla\phi-\frac{\phi}{\Psi}\nabla\Psi\right|^2\ge0
\end{align*}
proves the weighted inequality $pq\int_\Omega v^{p-1}\phi^2\le\int_\Omega u^{1-q}(\Delta\phi)^2$. For any $f\in C_c^\infty(\Omega)$, substitute $\phi=G_\Omega(\sqrt{qu^{q-1}}\,f)$. This yields
\begin{equation*}
\bigl\|\sqrt{pv^{p-1}}\,G_\Omega\sqrt{qu^{q-1}}\,f\bigr\|_2\le\|f\|_2\,.
\end{equation*}
Increasing the interior domains, positive Green convergence and Fatou prove~\eqref{eq:green-contraction} on nonnegative compactly supported inputs. Positivity and density then give the full $L^2$ contraction.

For the stated version of~\eqref{eq:Navier} on $D$, apply this contraction to $f=-q^{-1/2}u^{(1-q)/2}\Delta\phi$ and use $G_D(-\Delta\phi)=\phi$. Interior source truncations and uniqueness of the Dirichlet energy solution justify this identity in the weighted operator; the assumed right-hand side is exactly its required integrability. Thus no trace of a Jacobi supersolution on the boundary has been introduced.

For $\zeta\in C_c^\infty(D)$, the Picone calculation with $\zeta^2/\Phi$ and $\zeta^2/\Psi$ gives
\begin{align*}
\int_D|\nabla\zeta|^2
&\ge\frac12\int_D
\left(pv^{p-1}\frac{\Psi}{\Phi}+qu^{q-1}\frac{\Phi}{\Psi}\right)\zeta^2\\
&\ge\sqrt{pq}\int_D u^{(q-1)/2}v^{(p-1)/2}\zeta^2\,.
\end{align*}
This is the argument of~\cite[Lemma 3]{Cowan2013}, applied to the positive supersolutions on $D$. Density in $H_0^1$, followed by Fatou, proves~\eqref{eq:Cowan} for the stated zero-trace tests. All exhaustions take place inside $D$, so no boundary value of $\Phi$ or $\Psi$ is used.
\end{proof}

The following lemma characterizes stability on an open set in terms of Green operator estimates on its bounded interior subdomains. It will be used to establish the stability of the positive limit constructed in Section~\ref{sec:global}.

\begin{lemma}\label{lem:local-stability}
Let $n\ge1$, let $\Omega\subseteq\mathbb{R}_+^n$ be open, and let $a,b\in C^\infty(\Omega)$ be strictly positive. There exist strictly positive $\Phi,\Psi\in C^2(\Omega)$ with $-\Delta\Phi\ge a\Psi$ and $-\Delta\Psi\ge b\Phi$ in $\Omega$ if and only if
\[
\|\sqrt a\,G_D\sqrt b\|_{L^2(D)\to L^2(D)}\le1
\]
for every bounded smooth connected domain $D\Subset\Omega$.
\end{lemma}

\begin{proof}
The Picone argument in Lemma~\ref{lem:stability} applies to these positive coefficients and proves the forward implication. For the converse, we first construct positive pairs on bounded domains and then normalize them at one fixed interior point to take a limit.

On a bounded smooth connected domain $D\Subset\Omega$, the compact positive self-adjoint operator $K_D=\sqrt a\,G_D bG_D\sqrt a$ has a largest eigenvalue $\lambda_D>0$ and a strictly positive eigenfunction $f_D$. Indeed, its Rayleigh quotient does not decrease on replacing a function by its absolute value, and its positive kernel makes a nonnegative eigenfunction strictly positive. Set $\Psi_D=f_D/\sqrt a$ and $\Phi_D=G_D(a\Psi_D)$. The eigenfunction equation gives
\begin{align*}
    \begin{aligned}
        -\Delta\Phi_D&=a\Psi_D,&
        -\Delta\Psi_D&=b\Phi_D/\lambda_D &&\text{in }D,\\
        \Phi_D&=0,&\Psi_D&=0 &&\text{on }\partial D.
    \end{aligned}
\end{align*}
Elliptic regularity, applied successively to the two equations, gives $\Phi_D,\Psi_D\in C^2(\overline D)$. This construction applies for every $\lambda_D>0$. Since $K_D=(\sqrt a\,G_D\sqrt b)(\sqrt a\,G_D\sqrt b)^*$, the assumed contraction gives $\lambda_D\le1$.

Suppose first that $\Omega$ is connected. Choose bounded smooth connected domains $D_j$ with $D_j\Subset D_{j+1}\Subset\Omega$ and $\bigcup_jD_j=\Omega$, and fix $x_0\in D_1$. Write $\lambda_j=\lambda_{D_j}$. Testing the Rayleigh quotient of $K_{D_j}$ on the zero extension of an $L^2(D_1)$-normalized positive eigenfunction of $K_{D_1}$, Green domain monotonicity gives $\lambda_j\ge\lambda_1>0$. Thus $\lambda_1\le\lambda_j\le1$ for every $j$. Normalize the constructed pair $(\Phi_j,\Psi_j)$ so that $\Phi_j(x_0)+\Psi_j(x_0)=1$. Its sum $Z_j=\Phi_j+\Psi_j$ satisfies
\[
-\Delta Z_j=V_j Z_j,\qquad
0\le V_j=\frac{a\Psi_j+b\Phi_j/\lambda_j}{\Phi_j+\Psi_j}
\le a+b/\lambda_1.
\]
The right-hand bound is uniform on every compact subset of $\Omega$. The scalar Harnack inequality on fixed chains of interior balls, starting at $x_0$, therefore bounds $Z_j$, and hence both components, on each such compact set. Interior $W^{2,s}$ estimates with $s>n$, followed by Schauder estimates for the two equations, give a subsequence converging in $C^2_{\loc}(\Omega)$. Taking a further subsequence, let $\lambda_j\to\lambda_*\in[\lambda_1,1]$. The limits satisfy
\[
-\Delta\Phi=a\Psi,\qquad -\Delta\Psi=b\Phi/\lambda_*\ge b\Phi,
\qquad \Phi(x_0)+\Psi(x_0)=1.
\]
Both components are nonnegative and superharmonic. If either vanished at an interior point, the strong minimum principle would make it identically zero, and the positive coupling would force the other to vanish as well. The normalization excludes this possibility. Thus both components are strictly positive on $\Omega$.

If $\Omega$ is not connected, apply the construction separately on each connected component and combine the resulting functions. Regularity is local, and no boundary values are prescribed, so this gives $\Phi,\Psi\in C^2(\Omega)$ with the required inequalities. All estimates used above are interior estimates; $a$ and $b$ need not be bounded near $\partial\Omega$.
\end{proof}

For $n\ge2$ and $p,q\ge1$ with $pq>1$, Li and Souplet~\cite[Proposition 5.8 and Lemma 5.9]{LiSouplet2026} derive a weighted Bessel-resolvent estimate on $\RR^{n-1}$ from normal monotonicity, using sine averages of positive normal derivatives and a two-weight Schur test. The coefficients in their final estimate are determined by the normal limits. Here~\eqref{eq:green-contraction} concerns the Dirichlet Green operator on $D\subseteq\mathbb{R}_+^n$, with the original coefficients $pv^{p-1}$ and $qu^{q-1}$. The estimate follows from stability through a fourth-order Picone inequality, requires no normal monotonicity, and remains valid when an exponent is below one.

Lemma~\ref{lem:local-stability} also proves the converse: Green contraction bounds on bounded interior domains yield a single pair of positive supersolutions on the entire open set. This construction uses positive eigenfunctions and an exhaustion argument with Harnack's inequality. It applies to general smooth positive coefficients, which need not be bounded near the boundary, and gives an equivalence between stability and the interior Green bounds.

\section{From exterior stability to global stability}\label{sec:global}

We show that the existence of a positive solution that is stable outside a compact set implies the existence of a possibly different solution stable throughout the half-space. Starting from a solution that stable outside a compact set, we construct minimal solutions on expanding half-balls with boundary data taken from the original solution. Each minimal solution is stable, and a subsequence converges locally either to zero or to a positive solution.

If the limit is positive, the linearized coefficients converge locally uniformly, and the Green operator estimates pass to the limit on every bounded smooth connected domain $D\Subset\mathbb R_+^n$. Lemma~\ref{lem:local-stability} then gives stability on the whole half-space. We exclude the zero limit by showing that it would force the original solution to have finite nonlinear energy, contradicting a boundary-flux identity and Hopf's lemma. This reduction allows us to use global stability in Section~\ref{sec:operators}.

\begin{lemma}\label{lem:global-reduction}
Let $n\ge2$, $p\ge q>0$, and $pq>1$. If~\eqref{eq:LE} has a positive classical solution stable outside a compact subset of $\overline{\mathbb{R}_+^n}$, then it has a positive classical solution stable on all of $\mathbb{R}_+^n$.
\end{lemma}

\begin{proof}
Fix an exterior-stable solution $(u,v)$ and its compact exceptional set $K$. On every bounded smooth interior subdomain of $\mathbb{R}_+^n\setminus K$, Lemma~\ref{lem:stability} gives the Green contraction. Lemma~\ref{lem:local-stability} therefore supplies one positive Jacobi pair on the entire complement. Restricting this pair allows us to use Lemma~\ref{lem:stability} on the half-annuli below even when their closures meet the boundary. In $B_R^+$, let $H_R^u,H_R^v$ be the positive harmonic extensions of the curved-boundary data of $(u,v)$, with zero data on the flat boundary. Set $(U_{R,0},V_{R,0})=(H_R^u,H_R^v)$ and, for $k\ge0$,
\begin{equation}\label{eq:halfball-iteration}
 \begin{aligned}
 U_{R,k+1}&=H_R^u+G_{B_R^+}(V_{R,k}^p),\\
 V_{R,k+1}&=H_R^v+G_{B_R^+}(U_{R,k}^q).
 \end{aligned}
\end{equation}
The iterates increase below $(u,v)$ to a minimal positive solution $(U_R,V_R)$ with the same boundary values. Minimality here is among nonnegative supersolutions with these boundary values and bounded above by $(u,v)$, by positive Green comparison and induction in~\eqref{eq:halfball-iteration}. The construction uses only that both powers are increasing. It also gives~\eqref{eq:component-comparison} for $(U_R,V_R)$ by Lemma~\ref{lem:comparison}.

Each minimal pair is stable on every smooth $\Omega\Subset B_R^+$. We adapt the minimal-solution perturbation argument in~\cite[proof of Theorem 1.5, Section 5]{Montenegro2005}, constructing the perturbation from the compact Green operator in the proof of Lemma~\ref{lem:local-stability}. Take $a=pV_R^{p-1}$ and $b=qU_R^{q-1}$ on $\Omega$. If the largest eigenvalue of $K_\Omega$ were $\lambda>1$, that construction would give positive Dirichlet functions $\Phi,\Psi$ with $-\Delta\Phi=a\Psi$ and $-\Delta\Psi=b\Phi/\lambda$. Choose $c$ with $1/\lambda<c<1$. The pair $(U_R-\varepsilon c\Phi,V_R-\varepsilon\Psi)$ has residuals
\begin{align*}
    \begin{aligned}
        -\Delta(U_R-\varepsilon c\Phi)-(V_R-\varepsilon\Psi)^p
        &=\varepsilon(1-c)a\Psi+O(\varepsilon^2\Psi^2),\\
        -\Delta(V_R-\varepsilon\Psi)-(U_R-\varepsilon c\Phi)^q
        &=\varepsilon(c-1/\lambda)b\Phi+O(\varepsilon^2\Phi^2).
    \end{aligned}
\end{align*}
On $\overline\Omega$ both components have positive lower bounds. The coefficients $a,b$ are bounded below by positive constants, the eigenfunctions are bounded, and the Taylor remainders are uniform, including for $q<1$. Thus the residuals are nonnegative for sufficiently small $\varepsilon>0$. Extend the perturbation by the original minimal pair outside $\Omega$. Its interface measures have signs $-\varepsilon c\,\partial_\nu\Phi>0$ and $-\varepsilon\partial_\nu\Psi>0$, by Hopf's lemma on the smooth interior domain. The extended pair is therefore a distributional supersolution. It remains above the harmonic starting pair if $\varepsilon$ is decreased further. Iteration below this smaller supersolution contradicts minimality. Hence $\lambda\le1$ on every such $\Omega$. Lemma~\ref{lem:local-stability} then gives a positive Jacobi pair on all of $B_R^+$.

For each fixed $R$, the differences $u-U_R$ and $v-V_R$ solve Poisson equations with bounded right-hand sides and zero data on the entire half-ball boundary. Oddly extend these right-hand sides across $x_n=0$ and solve the Dirichlet problems on the full ball. Uniqueness and odd symmetry identify their restrictions with the two differences. The full-ball estimates therefore give $W^{2,s}$ regularity for every finite $s>1$ and bounded gradients on the closed half-ball. Moreover, $u/x_n$ and $v/x_n$ have positive minima on the curved hemisphere, including its equator by Hopf's lemma for the original solution on the flat boundary. Comparison with a sufficiently small multiple of the harmonic function $x_n$ gives
\[
c_R x_n\le U_R,V_R\le C_R x_n\qquad\text{in }B_R^+.
\]
No Hopf lemma at a corner is needed. The constants here may depend on $R$; they justify calculations on each fixed half-ball and are not used as uniform exhaustion estimates.

The original pair bounds the sources on every compact set. Interior Poisson estimates, followed by H\"older estimates for positive powers, therefore give subsequential convergence of $(U_R,V_R)$ in $C^2_{\loc}(\mathbb{R}_+^n)$ as $R\to\infty$. Every limit is either zero in both components or positive in both components. It has continuous zero boundary values by $0\le(U_R,V_R)\le(u,v)$. A positive limit is stable on all of $\mathbb{R}_+^n$ by the limiting argument following Lemma~\ref{lem:local-stability}. It remains to exclude a sequence $R_j\to\infty$ for which $(U_j,V_j):=(U_{R_j},V_{R_j})\to(0,0)$ locally in $\mathbb{R}_+^n$.

\emph{Step 1: Finite energy when $q\ge1$.} Choose $r$ so large that the exceptional compact set lies in $B_r$. In $D_j=B_{R_j}^+\setminus\overline B_r$. Let
\[
W=u-U_j,\quad Z=v-V_j,\quad f=u^q-U_j^q,\quad g=v^p-V_j^p.
\]
Write $Z=G_{D_j}f+h_Z$ and $W=G_{D_j}g+h_W$. The nonnegative harmonic remainders have data only on the fixed inner hemisphere. Comparison with the positive harmonic function $x_n|x|^{-n}$ gives
\begin{align}\label{eq:inner-remainders}
    0\le h_W+h_Z\le C x_n|x|^{-n}\qquad\text{in }D_j,
\end{align}
with $C$ independent of $j$. On the fixed inner hemisphere the boundary values satisfy $h_W+h_Z=W+Z\le u+v\le C x_n$, by local $C^1$ regularity of the original solution at the flat boundary and compactness. A suitable multiple of $x_n|x|^{-n}$ therefore dominates these data, without an appeal to Hopf's lemma at the equator.

Let $z_0=G_{D_j}f$. Applying~\eqref{eq:green-contraction} on $D_j$ to $(qu^{q-1})^{-1/2}f$, and then using the convexity of $s^q$ and symmetry of the Green kernel, gives
\begin{align*}
    \begin{aligned}
        \int_{D_j}pv^{p-1}z_0^2
        &\le\int_{D_j}\frac{f^2}{qu^{q-1}}
        \le\int_{D_j}fW\\
        &=\int_{D_j}g z_0+\int_{D_j}f h_W.
    \end{aligned}
\end{align*}
All these integrals are finite on the fixed half-annulus; in particular $f^2/(qu^{q-1})\le u^{q+1}/q$. Since $0\le g\le pv^{p-1}Z$, expansion of $Z=z_0+h_Z$ gives
\[
0\le\int_{D_j}(pv^{p-1}Z^2-gZ)
\le\int_{D_j}u^q h_W+2p\int_{D_j}v^p h_Z\le C.
\]
For the final bound, combine~\eqref{eq:inner-remainders} with~\eqref{eq:source-moments} on dyadic half-annuli. Their contributions are bounded by $CL^{-1-\beta}$ and $CL^{-1-\alpha}$, respectively. Fatou along the assumed zero-limit sequence gives $(p-1)\int_{\mathbb{R}_+^n\setminus B_r}v^{p+1}\le C$. Thus $v^{p+1}\in L^1(\mathbb{R}_+^n)$.

\emph{Step 2: Finite energy when $0<q<1$.} Fix a smooth radial cutoff $\chi$ equal to zero on $B_{2r}$ and one outside $B_{3r}$. Set
\[
\rho_j=U_j/u,\quad \sigma_j=V_j/v,\quad F_j=1-\rho_j^q,\quad
\widetilde W_j=uF_j/q,\quad Z_j=v-V_j.
\]
The function $(x,y)\mapsto x^{-(1-q)/q}y^p$ is convex on the positive quadrant because $p-1-(1-q)/q=p-1/q>0$. Its tangent inequality at $(1,1)$ yields
\[
\mathcal E_j=\rho_j^{q-1}\sigma_j^p-1-\frac{1-q}{q}F_j+p(1-\sigma_j)\ge0.
\]
A direct differentiation of the quotient $U_j/u$ gives
\begin{align}\label{eq:inverse-power-difference}
    \begin{aligned}
        -\Delta Z_j&=qu^{q-1}\widetilde W_j,\\
        -\Delta\widetilde W_j&=pv^{p-1}Z_j-\mathcal K_j,\\
        \mathcal K_j&=v^p\mathcal E_j+
        \frac{1-q}{q^2}u\frac{|\nabla F_j|^2}{1-F_j}\ge0.
    \end{aligned}
\end{align}
This change of variable retains the coefficient $qu^{q-1}$ of the original linearization.

The cutoff errors are uniformly bounded on the fixed inner half-annulus. To see this at the flat boundary, Lemma~\ref{lem:comparison} gives $V_j^p/U_j\le C U_j^{(pq-1)/(p+1)}\le C$ there. Interior Harnack and gradient estimates on balls of radius comparable to $x_n$ imply $|\nabla U_j|\le CU_j/x_n$. Therefore
\[
u^{1-q}U_j^{q-1}|\nabla U_j|\le Cu/x_n\le C.
\]
The remaining derivatives of $\widetilde W_j=(u-u^{1-q}U_j^q)/q$ are bounded for the same reason. Boundary Poisson estimates also bound $Z_j$ and $\nabla Z_j$ uniformly on this fixed half-annulus.

Apply~\eqref{eq:Navier} on $\widehat D_j:=B_{R_j}^+\setminus\overline B_{3r/2}$ to $\chi Z_j$. Since $q<1$, the weight $u^{1-q}$ is bounded on each fixed half-ball. Expanding only the derivatives of $\chi$ gives, with $f_j=u^q-U_j^q$,
\[
p\int_{\widehat D_j}\chi^2v^{p-1}Z_j^2\le\int_{\widehat D_j}\chi^2\widetilde W_j f_j+C.
\]
Cross-testing the two equations in~\eqref{eq:inverse-power-difference} gives
\begin{align*}
    \begin{aligned}
        \int_{\widehat D_j}\chi^2\widetilde W_j f_j
        &=p\int_{\widehat D_j}\chi^2v^{p-1}Z_j^2-\int_{\widehat D_j}\chi^2\mathcal K_jZ_j\\
        &\quad+2\int_{\widehat D_j}\chi(\widetilde W_j\nabla Z_j-Z_j\nabla\widetilde W_j)\cdot\nabla\chi.
    \end{aligned}
\end{align*}
Here and in the preceding display the domain is $\widehat D_j$. Both differences vanish on the outer hemisphere and on the boundary. The fixed-half-ball estimates above give $c_j x_n\le U_j,V_j\le C_j x_n$ and bounded gradients for each $j$. They also give $\chi Z_j\in H^2\cap H_0^1$ on the testing half-annulus. Thus $\widetilde W_j,Z_j=O_j(x_n)$, their gradients are bounded, and $\mathcal K_j=O_j(x_n^p+x_n^{-1})$. In particular $\mathcal K_jZ_j$ is integrable, and the flat-boundary cross term tends to zero after truncating at $x_n=\varepsilon$ and letting $\varepsilon\downarrow0$. The fixed-$j$ constants are used only for this justification. Uniform estimates are needed only on the support of $\nabla\chi$, where they were proved above. The boundary limit is taken for each fixed $j$ before $j\to\infty$. Cancellation of the two finite principal terms now yields $\int_{\widehat D_j}\chi^2\mathcal K_jZ_j\le C$.

Since $\rho_j,\sigma_j\to0$ locally, the nonnegative term $\rho_j^{q-1}\sigma_j^p$ can be retained without taking its limit, and
\[
\liminf_j\mathcal E_j\ge p-1/q>0.
\]
Fatou again gives $v^{p+1}\in L^1(\mathbb{R}_+^n)$.

\emph{Step 3: The boundary excludes finite energy.} By Lemma~\ref{lem:source} and positive Tonelli,
\[
\int_{\mathbb{R}_+^n}u^{q+1}=\int_{\mathbb{R}_+^n}u^qG_{\mathbb{R}_+^n}(v^p)
=\int_{\mathbb{R}_+^n}v^pG_{\mathbb{R}_+^n}(u^q)=\int_{\mathbb{R}_+^n}v^{p+1}<\infty.
\]
The zero boundary values and local $C^1$ regularity imply that odd reflection across $x_n=0$ produces the Poisson equations on $\RR^n$ without a boundary measure. We justify the global second-derivative estimate, including the possible polynomial ambiguity. Let $\widetilde u$ be the odd extension of $u$. The Calder\'on--Zygmund multipliers give $F_{ij}\in L^{(p+1)/p}(\RR^n)$ such that $D_{ij}\widetilde u-F_{ij}$ is a harmonic polynomial $P_{ij}$: its Laplacian is zero, and both distributions are tempered. For $\varphi\in C_c^\infty(\RR^n)$, let $\varphi_R(x)=R^{-n}\varphi(x/R)$. H\"older's inequality gives
\begin{align*}
    \begin{aligned}
        |\langle D_{ij}\widetilde u,\varphi_R\rangle|
        &\le C_\varphi R^{-2-n/(q+1)}\|\widetilde u\|_{L^{q+1}},\\
        |\langle F_{ij},\varphi_R\rangle|
        &\le C_\varphi R^{-np/(p+1)}\|F_{ij}\|_{L^{(p+1)/p}}.
    \end{aligned}
\end{align*}
Both terms tend to zero. A nonzero polynomial cannot have this property for every $\varphi$, since one can choose $\varphi$ with a nonzero moment against its highest homogeneous part. Hence $P_{ij}=0$, and $D^2u\in L^{(p+1)/p}(\mathbb{R}_+^n)$. Interchanging the components and exponents gives $D^2v\in L^{(q+1)/q}(\mathbb{R}_+^n)$.

The Gagliardo--Nirenberg inequality gives $\nabla u\in L^{r_u}$ and $\nabla v\in L^{r_v}$, where
\[
\frac1{r_u}=\frac12\left(\frac p{p+1}+\frac1{q+1}\right),\qquad
\frac1{r_v}=\frac12\left(\frac q{q+1}+\frac1{p+1}\right),\qquad
\frac1{r_u}+\frac1{r_v}=1.
\]
For the boundary argument we use the normal-translation form of the mixed Rellich identity for elliptic systems; see~\cite[Corollary 2.1, equation (2.4), p.~128]{Mitidieri1993}. The scalar normal-translation and boundary-flux argument appears in~\cite[proof of Proposition 3]{DupaigneFarinaPetitt2023}. We give the mixed tensor calculation needed for the present system; its justification uses the finite reaction energy and conjugate gradient integrability established above. The tensor
\[
T=\nabla u\otimes\nabla v+\nabla v\otimes\nabla u
-\left(\nabla u\cdot\nabla v-\frac{u^{q+1}}{q+1}-\frac{v^{p+1}}{p+1}\right)I
\]
belongs to $L^1(\mathbb{R}_+^n)$. Direct differentiation gives $\operatorname{div}T=0$, while $Te_n=u_{x_n} v_{x_n} e_n$ on the boundary. For a nonnegative cutoff $\eta_R$ equal to one on $B_R$ and supported in $B_{2R}$, integrate on a smooth bounded domain whose boundary agrees with $x_n=0$ near the support of $\eta_R$ and whose remaining boundary is disjoint from that support. The divergence theorem gives
\begin{align*}
    \begin{aligned}
        \int_{\partial \mathbb{R}_+^n}\eta_Ru_{x_n} v_{x_n}&=-\int_{\mathbb{R}_+^n}T:\nabla(\eta_Re_n),\\
        \left|\int_{\mathbb{R}_+^n}T:\nabla(\eta_Re_n)\right|
        &\le\frac C R\int_{B_{2R}^+\setminus B_R^+}|T|\longrightarrow0.
    \end{aligned}
\end{align*}
Hopf's lemma makes the left side bounded below by a positive integral on any fixed nonempty boundary disk. This is a contradiction. The identity uses normal translation and has no exceptional scaling exponent.

The zero limit has been excluded. A positive subsequential limit of the minimal half-ball solutions therefore proves the lemma. Every estimate above is valid for $n=2$ as well: the half-space Green kernel is logarithmic, and the exterior harmonic barrier in~\eqref{eq:inner-remainders} is still $x_n|x|^{-n}$.
\end{proof}

\section{Minimal solutions on finite strips}\label{sec:strips}

Fix $n\ge2$, $p\ge q>0$, $pq>1$, and a positive classical solution $(u,v)$ of~\eqref{eq:LE}. For $S>0$, write $\Sigma_S=\RR^{n-1}\times(0,S)$ and $G_S=G_{\Sigma_S}$. We will construct an auxiliary solution on $\Sigma_S$ with boundary values $(0,0)$ at $x_n=0$ and $(u,v)(\cdot,S)$ at $x_n=S$, and obtain bounds uniform in the tangential variable $x'$. No stability assumption is used in this section.

For nonnegative continuous data $g$ on $\RR^{n-1}$, let $\Pi_Sg$ denote the minimal nonnegative harmonic extension with boundary values $0$ at $x_n=0$ and $g$ at $x_n=S$. It is defined by increasing limits of bounded harmonic extensions of compactly supported truncations of $g$. When this limit is finite, minimality means that $\Pi_Sg\le w$ for every nonnegative harmonic function $w$ on $\Sigma_S$ with the same boundary values.

Green comparison and Lemma~\ref{lem:source}\textup{(iii)} give
\begin{equation*}
0\le G_S(v^p)\le u,\qquad 0\le G_S(u^q)\le v\quad\text{in }\Sigma_S\,.
\end{equation*}
Hence the functions
\begin{equation*}
b_S:=u-G_S(v^p),\qquad c_S:=v-G_S(u^q)
\end{equation*}
are well defined and nonnegative. The equations imply $\Delta b_S=\Delta c_S=0$ in $\Sigma_S$, so $b_S,c_S$ are the harmonic parts of $u,v$ after subtracting the Green potentials of their sources.

Starting from $(M_0,N_0)=(b_S,c_S)$, define, for every integer $j\ge0$,
\begin{align}\label{eq:Picard}
\begin{aligned}
M_{j+1}&=b_S+G_S(N_j^p),\\
N_{j+1}&=c_S+G_S(M_j^q)\,.
\end{aligned}
\end{align}
The following lemma shows that these iterates increase and remain bounded above by $(u,v)$, and that their limit is a positive solution with the prescribed boundary values.

\begin{lemma}\label{lem:strip-construction}
For every $S>0$, the iterates defined by~\eqref{eq:Picard} increase pointwise to a positive classical solution $(M_S,N_S)$ of
\begin{equation*}
\begin{aligned}
-\Delta M_S&=N_S^p,\qquad -\Delta N_S=M_S^q
&&\text{in }\Sigma_S,\\
M_S&=N_S=0
&&\text{on }\{x_n=0\},\\
(M_S,N_S)&=(u,v)
&&\text{on }\{x_n=S\}.
\end{aligned}
\end{equation*}
This solution satisfies $M_S\le u$ and $N_S\le v$ in $\Sigma_S$ and has the following properties.
\begin{enumerate}
\item[\textup{(i)}] For every nonnegative classical supersolution $(W_1,W_2)$ of the same system in $\Sigma_S$ with the same boundary values and $W_1\le u$, $W_2\le v$, we have $M_S\le W_1$ and $N_S\le W_2$ in $\Sigma_S$.

\item[\textup{(ii)}] For every $x'\in\RR^{n-1}$ and $0<x_n<S/2$,
\begin{equation}\label{eq:minimal-monotonicity}
(M_S)_{x_n}(x',x_n)\ge0,\qquad (N_S)_{x_n}(x',x_n)\ge0\,.
\end{equation}

\item[\textup{(iii)}] Let $0<h<S$, and let $(W_1,W_2)$ be a positive classical supersolution of the same system in $\Sigma_h$, continuous up to both horizontal boundaries and zero on $x_n=0$. If $W_1(\cdot,h)\ge M_S(\cdot,h)$ and $W_2(\cdot,h)\ge N_S(\cdot,h)$ on $\RR^{n-1}$, then $M_S\le W_1$ and $N_S\le W_2$ in $\Sigma_h$.
\end{enumerate}
\end{lemma}

\begin{proof}
We first establish a comparison principle on strips. Let $h>0$, and let $z\ge0$ be subharmonic in $\Sigma_h$, continuous up to $x_n=0,h$, and zero there. Suppose $z\le u$ in $\Sigma_h$. Its extension by zero to $\RR^n$ is subharmonic, by the truncation argument in the proof of Lemma~\ref{lem:source}. The mean-value inequality and~\eqref{eq:horizontal-L1} give, for every fixed $x\in\Sigma_h$,
\begin{equation*}
z(x)\le CR^{-n}\int_{B_R(x)\cap\Sigma_h}u
\le ChR^{-1-\alpha}\longrightarrow0
\qquad\text{as }R\to\infty\,.
\end{equation*}
Thus $z=0$. If $z\le v$, the last bound is $ChR^{-1-\beta}$, so the same conclusion holds.

\emph{Existence and \textup{(i)}.}
The functions $b_S,c_S$ are nonnegative and harmonic in $\Sigma_S$. Since $b_S\le u$ and $c_S\le v$, they extend continuously by zero to $x_n=0$. We next verify their values at $x_n=S$. Fix $\xi\in\RR^{n-1}$ and $r\in(0,S/4)$, and split each source into its restriction to $B_{2r}((\xi,S))\cap\Sigma_S$ and the complementary part. The local source is bounded. Comparison with a multiple of $x_n(S-x_n)/2$ shows that its Green potential tends to zero as $x_n\uparrow S$. For the complementary part, the kernel is positive and harmonic near $(\xi,S)$ and vanishes on $x_n=S$. Boundary Harnack estimates therefore give
\begin{equation*}
G_S(x,y)\le C(S-x_n)G_S(x_*,y)\,,
\end{equation*}
for $x=(x',x_n)\in B_r((\xi,S))\cap\Sigma_S$ and $y\in\Sigma_S\setminus B_{2r}((\xi,S))$, where $x_*=(\xi,S-r)$ and $C$ is independent of $y$. The inequalities $G_S(v^p)(x_*)\le u(x_*)$ and $G_S(u^q)(x_*)\le v(x_*)$ then show that the complementary potentials also tend to zero. Hence $b_S(\cdot,S)=u(\cdot,S)$ and $c_S(\cdot,S)=v(\cdot,S)$.

By the definition of the minimal harmonic extension,
\begin{equation*}
0\le\Pi_S[u(\cdot,S)]\le b_S\,.
\end{equation*}
The compactly supported approximations defining $\Pi_S$ give the lower boundary limit $u(\xi,S)$ at each point $(\xi,S)$, while $b_S$ gives the matching upper limit. At $x_n=0$, the extension tends to zero because it is bounded by $u$. Thus $b_S-\Pi_S[u(\cdot,S)]$ is nonnegative, harmonic, bounded by $u$, and zero on both horizontal boundaries. The comparison principle just proved makes this difference zero. Applying the same argument with $v$ gives
\begin{equation*}
b_S=\Pi_S[u(\cdot,S)],\qquad c_S=\Pi_S[v(\cdot,S)]\,.
\end{equation*}
In particular, $b_S,c_S>0$ in $\Sigma_S$ by the strong minimum principle.

Positivity of $G_S$ and monotonicity of the powers give, by induction in~\eqref{eq:Picard},
\begin{equation*}
0<M_j\le M_{j+1}\le u,\qquad
0<N_j\le N_{j+1}\le v
\qquad(j\ge0)\,.
\end{equation*}
The pointwise limits are therefore finite and positive. Monotone convergence gives
\begin{equation*}
M_S=b_S+G_S(N_S^p),\qquad
N_S=c_S+G_S(M_S^q)\,.
\end{equation*}
On each compact subset of $\Sigma_S$, the iterates are bounded above by $(u,v)$ and below by $(b_S,c_S)>0$. Interior elliptic estimates give local $C^2$ convergence, so the limits solve the system classically. The bounds $b_S\le M_S\le u$ and $c_S\le N_S\le v$ give the prescribed continuous boundary values.

Let $(W_1,W_2)$ satisfy the assumptions in~\textup{(i)}. Comparison on bounded cylinders, with zero lateral boundary data, followed by monotone exhaustion gives
\begin{equation*}
W_1\ge b_S+G_S(W_2^p),\qquad
W_2\ge c_S+G_S(W_1^q)\,.
\end{equation*}
Starting from $M_0=b_S\le W_1$ and $N_0=c_S\le W_2$, induction in~\eqref{eq:Picard} yields $M_j\le W_1$ and $N_j\le W_2$. Passing to the limit proves~\textup{(i)}.

\emph{Proof of \textup{(ii)}.}
Fix $0<\lambda\le S/2$. For a function $F$ on $\Sigma_S$, write $F^\lambda(x',x_n)=F(x',2\lambda-x_n)$ on $\Sigma_\lambda$. The differences $b_S^\lambda-b_S$ and $c_S^\lambda-c_S$ are harmonic, nonnegative at $x_n=0$, and zero at $x_n=\lambda$. Their negative parts are bounded by $u$ and $v$, respectively, so the comparison principle gives $b_S^\lambda\ge b_S$ and $c_S^\lambda\ge c_S$.

Suppose $M_j^\lambda\ge M_j$ and $N_j^\lambda\ge N_j$. The equations for the next iterates give
\begin{align*}
-\Delta(M_{j+1}^\lambda-M_{j+1})
&=(N_j^\lambda)^p-N_j^p\ge0,\\
-\Delta(N_{j+1}^\lambda-N_{j+1})
&=(M_j^\lambda)^q-M_j^q\ge0\,.
\end{align*}
These differences are nonnegative on the two horizontal boundaries of $\Sigma_\lambda$. Their negative parts are subharmonic and bounded by $M_{j+1}\le u$ and $N_{j+1}\le v$. The comparison principle therefore proves the reflection inequalities for $j+1$. Induction and passage to the limit give $M_S^\lambda\ge M_S$ and $N_S^\lambda\ge N_S$. For $0<t_1<t_2<S/2$, take $\lambda=(t_1+t_2)/2$ to obtain monotonicity in $x_n$. This proves~\eqref{eq:minimal-monotonicity}.

\emph{Proof of \textup{(iii)}.}
We first record the restriction identity
\begin{equation}\label{eq:strip-restriction}
Z|_{\Sigma_h}
=\Pi_h[Z(\cdot,h)]+G_h(f\mathbf1_{\Sigma_h}),
\qquad 0<h<S\,,
\end{equation}
for a finite function $Z=\Pi_Sg+G_Sf$ with nonnegative boundary data $g$ and source $f$. For bounded compactly supported data, both sides are bounded solutions of the same Poisson equation on $\Sigma_h$ with the same horizontal boundary values. Indeed, the Poisson extension is bounded by the supremum of its boundary data, and
\begin{equation*}
0\le G_Sf(x',x_n)\le\|f\|_\infty\,\frac{x_n(S-x_n)}2\,.
\end{equation*}
The positive and negative parts of the difference, extended by zero, are bounded subharmonic functions supported in $\overline{\Sigma_h}$. Their ball averages are $O(h/R)$ as $R\to\infty$, so both vanish. Increasing bounded compactly supported truncations of $g,f$ and applying monotone convergence to the Green and Poisson kernels proves~\eqref{eq:strip-restriction} for the data used here.

Now let $(W_1,W_2)$ satisfy the assumptions in~\textup{(iii)}. Comparison on bounded cylinders and monotone exhaustion give
\begin{align*}
W_1&\ge\Pi_h[M_S(\cdot,h)]+G_h(W_2^p),\\
W_2&\ge\Pi_h[N_S(\cdot,h)]+G_h(W_1^q)\,.
\end{align*}
Apply~\eqref{eq:strip-restriction} first to $M_0=b_S$ and $N_0=c_S$. Their values at $x_n=h$ are bounded by $M_S(\cdot,h)$ and $N_S(\cdot,h)$, so $M_0\le W_1$ and $N_0\le W_2$ on $\Sigma_h$. If these inequalities hold for the $j$th iterates, then~\eqref{eq:strip-restriction} gives
\begin{align*}
M_{j+1}
&=\Pi_h[M_{j+1}(\cdot,h)]+G_h(N_j^p\mathbf1_{\Sigma_h})\\
&\le\Pi_h[M_S(\cdot,h)]+G_h(W_2^p)\le W_1
\end{align*}
on $\Sigma_h$, and likewise $N_{j+1}\le W_2$. Induction and passage to the limit prove~\textup{(iii)}. This comparison uses bounded-domain maximum principles and monotone limits, so it applies even when the supersolutions do not belong to $H^1$ near $x_n=0$.
\end{proof}

Unlike Li and Souplet~\cite[Proposition 3.1 and Lemma 5.9]{LiSouplet2026}, who average affine interpolation errors and positive normal derivatives to obtain linear inequalities and a two-weight Schur contraction, we work directly with the minimal strip solution $(M_S,N_S)$, without assuming normal monotonicity of the original solution. Weighted averages of $M_S^{\min(q,1)}$ and $N_S$ yield nonlinear differential inequalities, which lead to a supersolution controlling the Picard iterates. The construction also applies to $0<q<1$, using a chain-rule gradient term and H\"older's inequality.

\begin{lemma}\label{lem:weighted-trace}
Let $(M_S,N_S)$ be given by Lemma~\ref{lem:strip-construction}, and set $\rho=\min(q,1)$. There are $c_0,C>0$, depending only on $n,p,q$, such that this pair satisfies
\begin{align}\label{eq:small-traces}
    \begin{aligned}
        \int_{B'_r(\xi)}M_S(x',x_n)\md x'&\le C(x_n/S)r^{n-1-\alpha},\\
        \int_{B'_r(\xi)}N_S(x',x_n)\md x'&\le C(x_n/S)^\rho r^{n-1-\beta}
    \end{aligned}
\end{align}
for all $S,r>0$, all $\xi\in\RR^{n-1}$, and $0<x_n\le c_0S$.
\end{lemma}

\begin{proof}
Let $h=S/[4(1+4/\sqrt\rho)]$ and $L=4h/\sqrt\rho$, so that $h+L=S/4$ and $16/L^2=\rho/h^2$. For $z=(x_n-h)/L$, define the probability density
\[
w(x_n)=\frac{30}{L}z^2(1-z)^2\qquad(h\le x_n\le h+L).
\]
Its integral is one, and $w=w'=0$ at both endpoints. A direct calculation gives
\[
w''+\rho h^{-2} w=\frac{30}{L^3}\bigl(4(z-1/2)^2+16(z-1/2)^4\bigr)\ge0.
\]
With $\md\nu=w(x_n)\md x_n$, set
\[
A(x')=\left(\int M_S(x',x_n)^\rho\md\nu(x_n)\right)^{1/\rho},\qquad
B(x')=\int N_S(x',x_n)\md\nu(x_n).
\]
All averages are over positive interior heights, so the functions are locally smooth and positive. In this proof $\mathbb E$ denotes this average. Since
\[
-\Delta M_S^\rho=\rho M_S^{\rho-1}N_S^p+\rho(1-\rho)M_S^{\rho-2}|\nabla M_S|^2,
\]
two integrations by parts in $x_n$, with no endpoint term, give
\begin{align}\label{eq:weighted-power}
    \begin{aligned}
        (-\Delta'+\rho h^{-2})\mathbb E M_S^\rho
        &=\rho\mathbb E(M_S^{\rho-1}N_S^p)
        +\rho(1-\rho)\mathbb E(M_S^{\rho-2}|\nabla M_S|^2)\\
        &\quad+\int_h^{h+L}(w''+\rho h^{-2} w)M_S^\rho\md x_n.
    \end{aligned}
\end{align}
The last term is nonnegative. The chain rule for $A=(\mathbb E M_S^\rho)^{1/\rho}$ creates a negative tangential gradient term. It is bounded by the tangential part of the positive gradient term in~\eqref{eq:weighted-power}, because
\[
|\nabla'\mathbb E M_S^\rho|^2\le\rho^2(\mathbb E M_S^\rho)
\mathbb E(M_S^{\rho-2}|\nabla'M_S|^2).
\]
It follows that $(-\Delta'+h^{-2})A\ge A^{1-\rho}\mathbb E(M_S^{\rho-1}N_S^p)$.

Here $p>1$ and $p\rho>1$. If $\rho<1$, then $(1-\rho)/(p-1)<\rho$. H\"older's inequality and the order of power means give
\begin{align*}
    \begin{aligned}
        B^p&\le\mathbb E(M_S^{\rho-1}N_S^p)(\mathbb E M_S^{(1-\rho)/(p-1)})^{p-1}\\
        &\le A^{1-\rho}\mathbb E(M_S^{\rho-1}N_S^p).
    \end{aligned}
\end{align*}
For $\rho=1$ this is Jensen's inequality. A linear average of the second equation with the same weight gives
\begin{align}\label{eq:averaged-supersolution}
    \begin{aligned}
        (-\Delta'+h^{-2})A&\ge B^p,\\
        (-\Delta'+\rho h^{-2})B&\ge\mathbb E M_S^q\ge A^q.
    \end{aligned}
\end{align}
The last step uses $q\ge\rho$. In particular, it does not reverse Jensen's inequality for a concave power.

For $0\le x_n\le h$, take the explicit quadratic barrier $k(x_n)=2x_n/h-x_n^2/h^2$. Then $0\le k\le1$, $k(h)=1$, and $-k''=2/h^2\ge h^{-2}k$. Consequently, for $0<x_n<h$,~\eqref{eq:averaged-supersolution} gives
\begin{align*}
    \begin{aligned}
        -\Delta(kA)&\ge kB^p\ge(k^\rho B)^p,\\
        -(k^\rho)''&=\frac{2\rho}{h^2}k^{\rho-1}+\rho(1-\rho)k^{\rho-2}(k')^2\ge\rho h^{-2} k^\rho,\\
        -\Delta(k^\rho B)&\ge k^\rho A^q\ge(kA)^q.
    \end{aligned}
\end{align*}
The averaging interval ends at $S/4$, where~\eqref{eq:minimal-monotonicity} applies. Thus $A\ge M_S(\cdot,h)$ and $B\ge N_S(\cdot,h)$. The comparison following~\eqref{eq:strip-restriction} proves $M_S\le kA$ and $N_S\le k^\rho B$ on $\Sigma_h$. These supersolutions are continuous and zero at the bottom. When $\rho\le1/2$, the second one need not be in $H^1$ there; only the positive Green comparison is being used.

Since $\rho\le1$, we have $A\le\mathbb E M_S\le\mathbb E u$ and $B\le\mathbb E v$. Apply~\eqref{eq:horizontal-L1} at every averaging height and use $k(x_n)\le C_q x_n/S$. This proves~\eqref{eq:small-traces} with $c_0=[4(1+4/\sqrt\rho)]^{-1}$.
\end{proof}

\begin{lemma}\label{lem:strip-decay}
For the pair in Lemma~\ref{lem:strip-construction}, there exist $c,C>0$, depending only on $n,p,q$, such that for every $S>0$,
\begin{align}\label{eq:strip-decay}
    \begin{aligned}
        M_S(x',x_n)&\le C x_n S^{-\alpha-1},\\
        N_S(x',x_n)&\le C x_n S^{-\beta-1}
    \end{aligned}
    \qquad(x'\in\RR^{n-1},\ 0<x_n\le cS).
\end{align}
\end{lemma}

\begin{proof}
To convert Lemma~\ref{lem:weighted-trace} into a pointwise bound, we use the local doubling selection of~\cite[Lemma 5.1]{PolacikQuittnerSouplet2007}, with a finite-ball argument in place of a limiting Liouville theorem. Let $Q=M_S^{1/\alpha}+N_S^{1/\beta}$. Suppose that $x_0=(\xi_0,t_0)$, $t_0\le c_0S/4$, and $t_0Q(x_0)>64$. Repeatedly move to a point within distance $2/Q$ where $Q$ more than doubles, if such a point exists. The sum of all displacements is less than $4/Q(x_0)$, so the points and comparison balls remain in an interior compact neighborhood of $x_0$. Continuity forces the procedure to stop at a point $z$ with $Q(z)\ge Q(x_0)$ and $Q\le2Q(z)$ on $B_{2/Q(z)}(z)$.

With $\lambda=Q(z)^{-1}$, the rescaled pair $U(y)=\lambda^\alpha M_S(z+\lambda y)$, $V(y)=\lambda^\beta N_S(z+\lambda y)$ solves the original equations on $B_2$, is bounded there by $2^\alpha,2^\beta$, and satisfies $U(0)^{1/\alpha}+V(0)^{1/\beta}=1$. Bounded sources and interior Poisson estimates give a uniform gradient bound on $B_1$. Hence there are fixed $\delta\in(0,1/2)$ and $c_1>0$ such that
\begin{align}\label{eq:doubling-mass}
    \int_{B_\delta}(U+V)\ge c_1.
\end{align}
This uses only boundedness of the positive powers, not differentiability of $s^q$ at zero.

All horizontal slices of $B_{\delta\lambda}(z)$ have height between $t_0/2$ and $3t_0/2$. Apply~\eqref{eq:small-traces} with horizontal radius $\delta\lambda$ and integrate over its vertical extent. Since $x_n/S\le(x_n/S)^\rho$ here, scaling gives
\[
\int_{B_\delta}U\le C(t_0/S)^\rho\delta^{n-\alpha},\qquad
\int_{B_\delta}V\le C(t_0/S)^\rho\delta^{n-\beta}.
\]
Choose $0<c_2\le\min(c_0/4,1/4)$ such that $C c_2^\rho(\delta^{n-\alpha}+\delta^{n-\beta})<c_1$. For $t_0\le c_2S$ the last two bounds then contradict~\eqref{eq:doubling-mass}. Therefore $x_n Q(x',x_n)\le64$ for $x_n\le c_2S$. Apply this estimate at height $c_2S$, then use~\eqref{eq:minimal-monotonicity}, to obtain $M_S\le CS^{-\alpha}$ and $N_S\le CS^{-\beta}$ throughout $0<x_n\le c_2S$.

In each half-ball centered on the bottom with radius $c_2S/4$, the sources are bounded by $CS^{-\alpha-2}$ and $CS^{-\beta-2}$. The zero-boundary Poisson gradient estimate on the concentric half-ball of radius $c_2S/8$ gives $|\nabla M_S|\le CS^{-\alpha-1}$ and $|\nabla N_S|\le CS^{-\beta-1}$. Integrating from the zero bottom value proves~\eqref{eq:strip-decay}, after decreasing $c$. The trace estimate was used at arbitrarily small horizontal radii, not only at radii comparable to the height.
\end{proof}

\section{Energy estimates and the finite-height argument}\label{sec:operators}\label{sec:iteration}

For the rest of the proof, fix a positive classical solution $(u,v)$ of~\eqref{eq:LE} that is stable throughout $\mathbb{R}_+^n$, with $n\ge2$, $p\ge q>0$, and $pq>1$. Lemma~\ref{lem:global-reduction} reduces the remaining nonexistence argument to this case. We first estimate the integrals of $u^{q+1}$ and $v^{p+1}$. These estimates give polynomial growth of the measure $pv^{p+1}\md x$ on cylinders $B'_R(0)\times(0,h)$ as $R\to\infty$, for each fixed $h>0$.

We next compare the $L^2$ norms and $L^1$ masses of finite sums of iterates of a positive self-adjoint Green contraction, with respect to $pv^{p+1}\md x$. The norm estimates, tangential localization, and polynomial measure growth give mass upper bounds. Strip decay and nonlinear comparison yield incompatible lower bounds after finitely many fixed height increases.

\begin{lemma}\label{lem:energy}
For every $R>0$,
\begin{align}\label{eq:energy-growth}
    \int_{B_R^+}(u^{q+1}+v^{p+1})\le C R^{n-\alpha-\beta-2},
\end{align}
where $C$ depends only on $n,p,q$.
\end{lemma}

\begin{proof}
\emph{Step 1: Both exponents at least one.} We adapt the single-component power test and component comparison in~\cite[proof of Lemma 4]{Cowan2013}. After interchanging the two components in that argument, we take its testing exponent to be $(q+1)/2$. Here Lemma~\ref{lem:comparison} supplies the comparison up to the flat boundary, including $q=1$. Related two-component tests appear in~\cite[proof of Lemma 3.16]{HuangZou2026} and~\cite[Lemma 2.6]{LiSouplet2026}; we use the component comparison before estimating the cutoff errors, so one test suffices.

Scaling reduces the estimate to $R=1$. For $0\le\eta\in C_c^\infty(B_2)$, Let
\[
E=\int_{\mathbb{R}_+^n}u^qv^p\eta^2,\qquad
M_u=\int_{\mathbb{R}_+^n}u^{q+1}\bigl(|\nabla\eta|^2+|\Delta(\eta^2)|\bigr).
\]
The functions $u^q\eta^2$ and $u^{(q+1)/2}\eta$ belong to $H_0^1(\mathbb{R}_+^n)$ and have compact support in $\overline{\mathbb{R}_+^n}$, by the local boundary regularity and $q\ge1$. Testing the first equation with $u^q\eta^2$ gives
\[
q\int_{\mathbb{R}_+^n}u^{q-1}|\nabla u|^2\eta^2
=E+\frac1{q+1}\int_{\mathbb{R}_+^n}u^{q+1}\Delta(\eta^2).
\]
Expanding the gradient and integrating its mixed term by parts then yields
\begin{align}\label{eq:power-test}
    \begin{aligned}
        \int_{\mathbb{R}_+^n}|\nabla(u^{(q+1)/2}\eta)|^2
        &=\frac{(q+1)^2}{4q}E+\frac{1-q}{4q}\int_{\mathbb{R}_+^n}u^{q+1}\Delta(\eta^2)\\
        &\quad+\int_{\mathbb{R}_+^n}u^{q+1}|\nabla\eta|^2.
    \end{aligned}
\end{align}
The flat-boundary terms vanish because $u=0$ and $u^{q+1}=O(x_n^{q+1})$ there. By~\eqref{eq:component-comparison}, $u^{(q+1)/2}\ge\sqrt{(q+1)/(p+1)}\,v^{(p+1)/2}$. Apply~\eqref{eq:Cowan} to $u^{(q+1)/2}\eta$ and use~\eqref{eq:power-test} to obtain
\begin{align*}
    \begin{aligned}
        \sqrt{\frac{pq(q+1)}{p+1}}\,E
        &\le\sqrt{pq}\int_{\mathbb{R}_+^n}u^{(3q+1)/2}v^{(p-1)/2}\eta^2\\
        &\le\frac{(q+1)^2}{4q}E+CM_u.
    \end{aligned}
\end{align*}
The coefficient on the left is strictly larger than $(q+1)^2/(4q)$, since
\[
\frac{16pq^3}{(p+1)(q+1)^3}
\ge16\left(\frac q{q+1}\right)^4\ge1,
\]
and equality throughout would require $p=q=1$, contrary to $pq>1$. Hence $E\le CM_u$, with $C$ depending only on $p,q$. Formula~\eqref{eq:power-test} consequently bounds $\int_{\mathbb{R}_+^n}|\nabla(u^{(q+1)/2}\eta)|^2$ by $CM_u$.

Choose a fixed $\vartheta>1$, with $\vartheta\le n/(n-2)$ when $n>2$. The zero extension of $u^{(q+1)/2}\eta$ across the flat boundary belongs to $H^1$, so Sobolev's inequality gives
\[
\|u^{q+1}\eta^2\|_{L^\vartheta(\mathbb{R}_+^n)}\le C M_u.
\]
In dimension two any fixed finite $\vartheta>1$ is permitted. Take $\eta=\chi^L$, where $0\le\chi\le1$, $\chi=1$ on $B_1$, and $\chi\in C_c^\infty(B_2)$. Set $X=\int_{\mathbb{R}_+^n}u^{(q+1)\vartheta}\chi^{2L\vartheta}$. The $L^1$ bound in~\eqref{eq:source-moments} and H\"older's inequality imply
\[
X^{1/\vartheta}\le C X^{q/((q+1)\vartheta-1)},
\]
provided the integer $L$ satisfies $L\ge((q+1)\vartheta-1)/(\vartheta-1)$. Since $\vartheta q/((q+1)\vartheta-1)<1$, this bounds $X$ by a parameter-dependent constant. It bounds the unit-scale integral of $u^{q+1}$, and~\eqref{eq:component-comparison} then bounds that of $v^{p+1}$.

\emph{Step 2: A sublinear exponent.} Assume $0<q\le1$. The use of the fourth-order stability form and cutoff energy follows the approach of~\cite[Lemmas 2.1--2.3]{MtiriYe2019}; we give the calculation on the half-space, including the boundary terms. For a nonnegative smooth compact cutoff $\phi$ on $\overline{\mathbb{R}_+^n}$, multiplying $-\Delta v=u^q$ by $u^{1-q} v\phi$ and integrating the mixed gradient term once more gives
\begin{align*}
    \begin{aligned}
        &\int_{\mathbb{R}_+^n}u^{1-q}|\nabla v|^2\phi
        +\frac{q(1-q)}{2}\int_{\mathbb{R}_+^n}u^{-q-1}v^2|\nabla u|^2\phi
        +\frac{1-q}{2}\int_{\mathbb{R}_+^n}u^{-q}v^{p+2}\phi\\
        &\qquad=\int_{\mathbb{R}_+^n}uv\phi-\int_{\mathbb{R}_+^n}u^{1-q} v\nabla v\cdot\nabla\phi
        +\frac{1-q}{2}\int_{\mathbb{R}_+^n}u^{-q}v^2\nabla u\cdot\nabla\phi .
    \end{aligned}
\end{align*}
The boundary fluxes are $O(x_n^{2-q})$ on each compact boundary patch and tend to zero. Two applications of Young's inequality, absorbing the positive gradient terms on the left, yield
\begin{align}\label{eq:weighted-Caccioppoli}
    \int_{\mathbb{R}_+^n}u^{1-q}|\nabla v|^2\phi
    \le2\int_{\mathbb{R}_+^n}uv\phi+C_q\int_{\mathbb{R}_+^n}u^{1-q}v^2\frac{|\nabla\phi|^2}{\phi}.
\end{align}
The quotient is used only where $\phi>0$ and extended by zero in the cutoff constructions below. When $q=1$, all terms multiplied by $1-q$ disappear and the same estimate is the usual Caccioppoli inequality.

Choose $0\le\eta\le1$, equal to one on $B_R$ and supported in $B_{2R}$, with derivatives bounded at scales $R^{-1}$ and $R^{-2}$. Let $m=(q+1)/q$ be the conjugate exponent of $q+1$, choose an integer $k$ sufficiently large, and set
\[
E=\int_{\mathbb{R}_+^n}u^{q+1}\eta^{2k},\qquad J=\int_{\mathbb{R}_+^n}v^{p+1}\eta^{2k},\qquad
\mathcal R=R^{-2}\int_{\mathbb{R}_+^n}uv\eta^{2k-2}
+R^{-4}\int_{\mathbb{R}_+^n}u^{1-q}v^2\eta^{2k-4}.
\]
Apply~\eqref{eq:Navier} to $v\eta^k$. The expansion of its Laplacian, followed by~\eqref{eq:weighted-Caccioppoli} with $\phi=\eta^{2k-2}$, gives
\[
pJ\le\frac{1+\varepsilon}{q}E+C_{\varepsilon,k}\mathcal R.
\]
Testing the two equations against $v\eta^{2k}$ and $u\eta^{2k}$, respectively, and subtracting gives
\begin{align*}
    \begin{aligned}
        E-J
        &=\int_{\mathbb{R}_+^n}(u\nabla v-v\nabla u)\cdot\nabla(\eta^{2k})\\
        &=2\int_{\mathbb{R}_+^n}u\nabla v\cdot\nabla(\eta^{2k})
        +\int_{\mathbb{R}_+^n}uv\,\Delta(\eta^{2k}).
    \end{aligned}
\end{align*}
The second equality follows by one further integration by parts; both components vanish on the boundary. Young's inequality and~\eqref{eq:weighted-Caccioppoli}, with $\phi=\eta^{2k-2}$, therefore yield, for every $\delta>0$,
\[
|E-J|\le\delta E+C_{\delta,k}\mathcal R.
\]

Choose $\varepsilon,\delta>0$ such that $pq(1-\delta)>1+\varepsilon$. Combining the last two estimates yields $E+J\le C\mathcal R$.

Young's inequality now gives, for arbitrarily small $\nu>0$,
\[
\mathcal R\le\nu E+C_{\nu,k}R^{-2m}\int_{\mathbb{R}_+^n}v^m\eta^{2k-2m}.
\]
For example, the two powers of $u$ on the right before Young correspond to the fractions $1/(q+1)$ and $(1-q)/(q+1)$ of $E$. Both leave the power $m=(q+1)/q$ of $v$. At $q=1$ the second term is already the required remainder. Absorbing $\nu E$ and taking $k\ge m(p+1)/(p+1-m)$ gives
\[
E+J\le CR^{-2m}J^{m/(p+1)}R^{n(1-m/(p+1))}.
\]
Since $p+1>m$ and $2m/(p+1-m)=\beta$, this proves~\eqref{eq:energy-growth}. All tests have zero value on the flat boundary; their normal derivatives need not vanish. On compact boundary patches $u^{1-q}$ is bounded, and zero-trace $H^2$ approximation justifies~\eqref{eq:Navier} for $v\eta^k$. The displayed local boundary orders justify the integrations by parts. Finally, the scaling identities give the exponent in~\eqref{eq:energy-growth} in both cases.
\end{proof}

To compare different strip heights in one $L^2$ space, we use the fixed measure $\md\mu=pv^{p+1}\md x$. This weight makes the operators below self-adjoint, while~\eqref{eq:energy-growth} gives the polynomial measure bounds needed for the mass estimates. Composing the two linearized Green maps and writing a perturbation of $v$ as $vf$ leads to the operator
\begin{equation}\label{eq:Q-definition}
Q_hf=\frac{pq}{v}G_h\bigl(u^{q-1}G_h(v^pf)\bigr),\qquad h>0\,.
\end{equation}
We set $Q_hf=0$ outside $\Sigma_h$. For nonnegative measurable $f$, the Green integrals are interpreted by monotone convergence and may take the value $+\infty$. We use the kernel notation $Q_hf(x)=\int_{\mathbb{R}_+^n}f(y)Q_h(x,\md y)$.

For the rest of the paper, write $\langle f,g\rangle=\int_{\mathbb{R}_+^n}f\overline g\md\mu$ whenever the integral is absolutely convergent, and also allow the value $+\infty$ when $f,g\ge0$. The norms $\|\cdot\|_1$ and $\|\cdot\|_2$ denote the $L^1(\mu)$ and $L^2(\mu)$ norms, respectively. Operator norms without a specified space are taken on $L^2(\mathbb{R}_+^n,\mu)$.

The following properties will be used to compare operators at different strip heights and to estimate the norms of finite sums of iterates and their mass outside large horizontal cylinders.

\begin{lemma}\label{lem:operator}
The operators $Q_h$, $h>0$, have the following properties.
\begin{enumerate}
\item[\textup{(i)}] Each $Q_h$ has a nonnegative kernel symmetric with respect to $\mu$ and is a self-adjoint positive semidefinite contraction on $L^2(\mathbb{R}_+^n,\mu)$.

\item[\textup{(ii)}] For every $h>0$ and $x\in\mathbb{R}_+^n$,
\begin{equation*}
\int_{\mathbb{R}_+^n}Q_h(x,\md y)=Q_h1(x)\le pq\,.
\end{equation*}
For $0<h\le s$ and every nonnegative measurable function $f$,
\begin{equation*}
Q_hf\le Q_sf\,.
\end{equation*}

\item[\textup{(iii)}] There exist constants $a_n,C_n>0$, depending only on $n$, such that for every $h>0$ and $x\in\mathbb{R}_+^n$,
\begin{equation}\label{eq:exponential-row}
\int_{\mathbb{R}_+^n}
\exp\!\left(\frac{a_n|x'-y'|}{h}\right)Q_h(x,\md y)
\le C_n pq\,.
\end{equation}
\end{enumerate}
\end{lemma}

\begin{proof}
\emph{(i)}
Fix $h>0$ and set
\begin{equation*}
S_h=\sqrt{pv^{p-1}}\,G_h\sqrt{qu^{q-1}}\,.
\end{equation*}
By Lemma~\ref{lem:stability}, $S_h$ is a contraction on $L^2(\mathbb{R}_+^n,\md x)$, with all strip operators extended by zero outside $\Sigma_h$. Multiplication by $\sqrt p\,v^{(p+1)/2}$ maps $L^2(\mathbb{R}_+^n,\mu)$ unitarily onto $L^2(\mathbb{R}_+^n,\md x)$ and transforms $Q_h$ into $S_hS_h^*$. Hence $Q_h$ is self-adjoint, positive semidefinite, and has norm at most one.

The kernel is nonnegative, and direct substitution in~\eqref{eq:Q-definition} gives
\begin{equation*}
\mu(\md x)Q_h(x,\md y)
=p^2q\,v(x)^pv(y)^p
\left(\int_{\Sigma_h}G_h(x,z)u(z)^{q-1}G_h(z,y)\md z\right)
\md x\md y\,.
\end{equation*}
Since $G_h(x,z)=G_h(z,x)$, this measure is symmetric in $x,y$. This proves~\textup{(i)}.

\emph{(ii)}
Green comparison gives $G_h(v^p)\le u$ and $G_h(u^q)\le v$. Therefore,
\begin{equation*}
Q_h1
=\frac{pq}{v}G_h\bigl(u^{q-1}G_h(v^p)\bigr)
\le\frac{pq}{v}G_h(u^q)\le pq\,.
\end{equation*}
The total kernel mass equals $Q_h1$ by the kernel representation. For $0<h\le s$ and $f\ge0$, apply domain monotonicity to the inner and outer Green operators in~\eqref{eq:Q-definition}. Positivity of the coefficients then gives $Q_hf\le Q_sf$, including for nonnegative measurable $f$ by monotone convergence.

\emph{(iii)}
We first prove that, for constants $a_n,C_n>0$ depending only on $n$,
\begin{equation*}
\exp\!\left(\frac{a_n|x'-y'|}{h}\right)G_h(x,y)
\le C_nG_{\mathbb{R}_+^n}(x,y),
\qquad x,y\in\Sigma_h,\quad x\ne y\,.
\end{equation*}
Write $x=(x',x_n)$, $y=(y',s)$, and $\rho=|x'-y'|$. If $\rho\le2h$, this follows from domain monotonicity. Suppose $\rho>2h$. Let $k_h(\theta,x_n,s)$ denote the Dirichlet heat kernel of the interval $(0,h)$. The product formula for the strip heat kernel gives
\begin{equation*}
G_h(x,y)=\int_0^\infty
(4\pi\theta)^{-(n-1)/2}
\exp\!\left(-\frac{\rho^2}{4\theta}\right)
k_h(\theta,x_n,s)\md\theta\,.
\end{equation*}
Comparison with the half-line heat kernel for $\theta\le h^2$, and the sine-series expansion for $\theta\ge h^2$, give
\begin{equation*}
k_h(\theta,x_n,s)\le
\begin{cases}
C x_n s\,\theta^{-3/2},&0<\theta\le h^2,\\
C x_n s\,h^{-3}\exp(-\pi^2\theta/h^2),&\theta\ge h^2.
\end{cases}
\end{equation*}
Substitution in the preceding integral yields
\begin{equation*}
G_h(x,y)\le C_n x_n s\,h^{-n}\exp(-c_n\rho/h)\,.
\end{equation*}
The explicit half-space Green kernel, including its logarithmic form when $n=2$, satisfies
\begin{equation*}
G_{\mathbb{R}_+^n}(x,y)\ge c_n x_n s\,\rho^{-n}
\qquad(\rho>2h)\,.
\end{equation*}
Choose $0<a_n<c_n$. Since $(\rho/h)^n e^{-(c_n-a_n)\rho/h}$ is bounded for $\rho>2h$, the two estimates give the claimed kernel comparison.

For $x\in\Sigma_h$, use this comparison on both Green kernels in the representation of $Q_h$. The inequality $|x'-y'|\le|x'-z'|+|z'-y'|$ and Tonelli's theorem give
\begin{align*}
&\int_{\mathbb{R}_+^n}
\exp\!\left(\frac{a_n|x'-y'|}{h}\right)Q_h(x,\md y)\\
&\qquad\le\frac{C_n pq}{v(x)}
G_{\mathbb{R}_+^n}
\bigl(u^{q-1}G_{\mathbb{R}_+^n}(v^p)\bigr)(x)\\
&\qquad\le\frac{C_n pq}{v(x)}
G_{\mathbb{R}_+^n}(u^q)(x)
\le C_n pq\,,
\end{align*}
where we used $G_{\mathbb{R}_+^n}(v^p)\le u$ and $G_{\mathbb{R}_+^n}(u^q)\le v$. For $x\notin\Sigma_h$, the kernel is zero. This proves~\eqref{eq:exponential-row}.
\end{proof}

Let
\begin{equation}\label{eq:operator-ratios}
l_h=N_h/v,\qquad r_h=1-l_h\quad\text{on }\Sigma_h,\qquad \tau=x_n/v,\qquad\gamma=\min(1,1/q).
\end{equation}
Extend $l_h$ and $r_h$ by zero outside $\Sigma_h$; the identity $r_h=1-l_h$ is used only inside $\Sigma_h$.

\begin{lemma}\label{lem:nonlinear-operator}
For every $h>0$, there are $\delta_*>0$, depending only on $p,q$, and $K_h<\infty$, depending only on $p,q,h$, such that on $\Sigma_h$,
\begin{align}\label{eq:nonlinear-feedback}
    r_h+\delta_*Q_h(r_h^2)\le Q_hr_h,
\end{align}
\begin{align}\label{eq:helper}
    Q_hl_h\le pq l_h^{1/(pq)},\qquad
    Q_h1\ge pq(1-l_h^\gamma),
\end{align}
and
\begin{align}\label{eq:Foster}
    Q_h\tau\le K_h\tau^{1/(pq)}.
\end{align}
\end{lemma}

\begin{proof}
To prove~\eqref{eq:nonlinear-feedback}, write $b_h=u-G_h(v^p)$, $c_h=v-G_h(u^q)$, and
\[
\mathcal F(w)=c_h+G_h\bigl[(b_h+G_h(w^p))^q\bigr].
\]
Both $v$ and $N_h$ are fixed points. For each fixed spatial point $x\in\Sigma_h$, the inner expression is the $pq$-th power of an $L^p$ norm with measure $G_h(x,y)\md y$, after adjoining a unit-mass atom whose function value is $b_h(x)^{1/p}$. This observation explains why joint convexity holds even for $q<1$. For completeness, let $z=v-N_h$, $w_s=v-sz$, and $U_s=b_h+G_h(w_s^p)$. Differentiation gives
\begin{align*}
    \begin{aligned}
        \frac{\md^2}{\md s^2}U_s^q
        &=pq(p-1)U_s^{q-1}G_h(w_s^{p-2}z^2)\\
        &\quad+p^2q(q-1)U_s^{q-2}[G_h(w_s^{p-1}z)]^2.
    \end{aligned}
\end{align*}
When $q<1$, Cauchy--Schwarz bounds the square by $U_sG_h(w_s^{p-2}z^2)$ and leaves the positive coefficient $pq(pq-1)$. When $q\ge1$, the second term is nonnegative. On $0\le s\le1/2$, $w_s\ge v/2$ and $U_s\ge2^{-p}u$. Integration of the Taylor remainder therefore gives~\eqref{eq:nonlinear-feedback}, for example with
\[
\delta_*=\min\!\left\{\frac14,\,
\frac38\min(p-1,pq-1)\,2^{-p(q-1)_+-(p-2)_+}\right\}.
\]
The derivatives on this half-interval are dominated by a constant times $u^q$, whose Green potential is finite. Convexity extends the half-interval remainder bound to the entire difference. Thus no differentiation at a possibly zero endpoint is required.

For the first inequality in~\eqref{eq:helper}, two applications of H\"older give
\begin{align*}
    \begin{aligned}
        G_h(v^{p-1}N_h)&\le u^{1-1/p}M_h^{1/p},\\
        G_h(u^{q-1/p}M_h^{1/p})&\le v^{1-1/(pq)}N_h^{1/(pq)}.
    \end{aligned}
\end{align*}
The outer exponents are legitimate because $q-1/p>0$. The deficit in total kernel mass is exactly
\[
1-Q_h1/(pq)=\frac{c_h+G_h(u^{q-1}b_h)}v.
\]
If $q\le1$, use $u^{q-1}b_h\le M_h^{q-1}b_h\le M_h^q$ to bound it by $l_h$. If $q>1$, form the probability measure with an atom of mass $c_h/v$ and the density $G_h(x,y)u(y)^q/v(x)$. Give the atom value one and the other points value $b_h/u$. Its $q$-th moment is at most $(c_h+G_h(M_h^q))/v=l_h$. Jensen's inequality bounds its first moment by $l_h^{1/q}$, as required.

Finally, the horizontal integral of the strip Green kernel is the one-dimensional interval kernel. Hence, for $a>0$,
\[
G_h(x_n^a)=\frac{h^{a+1}x_n-x_n^{a+2}}{(a+1)(a+2)}
\le\frac{h^{a+1}x_n}{(a+1)(a+2)}.
\]
Two applications of H\"older give
\begin{align*}
    \begin{aligned}
        G_h(v^{p-1}x_n)&\le u^{1-1/p}[G_h(x_n^p)]^{1/p},\\
        G_h(u^{q-1/p}x_n^{1/p})&\le v^{1-1/(pq)}[G_h(x_n^q)]^{1/(pq)}.
    \end{aligned}
\end{align*}
Substitution into~\eqref{eq:Q-definition} proves~\eqref{eq:Foster}, with
\[
K_h=pq[(p+1)(p+2)]^{-1/p}[(q+1)(q+2)]^{-1/(pq)}
h^{(p+1)/p+(q+1)/(pq)}.
\]
All inequalities involve finite positive potentials. The possibly singular function $\tau$ is never used as a differential test function.
\end{proof}

To raise the strip height repeatedly, we need bounds on $\sum_{k\ge0}\langle f,Q_h^kf\rangle$ for the new nonnegative sources $f$. Fix a nonzero nonnegative $g_0\in C_c^\infty(\mathbb{R}_+^n)$ and use it at every height. For every $h>0$ with $\supp g_0\subseteq\Sigma_h$, define
\begin{align*}
F_h&=\sum_{k=0}^\infty Q_h^kg_0,\\
E(h)&=\langle g_0,F_h\rangle
=\sum_{k=0}^\infty\langle g_0,Q_h^kg_0\rangle\,.
\end{align*}
All series are interpreted as monotone limits. Although $F_h$ need not belong to $L^2(\mu)$, the next lemma shows that $E(h)$ is finite and that $E(s)-E(h)$ bounds the quadratic series of the sources produced when the height increases from $h$ to $s$.

\begin{lemma}\label{lem:source-energy}
For every $h>0$ with $\supp g_0\subseteq\Sigma_h$, we have $E(h)<\infty$ and $F_h\in L^1_{\loc}(\Sigma_h)$. For $s>h$, set $B_{h,s}=(Q_s-Q_h)F_h$, where the difference acts through the nonnegative kernel $Q_s(x,\md y)-Q_h(x,\md y)$. Then
\begin{equation}\label{eq:resolvent-difference}
\sum_{k=0}^\infty Q_s^kB_{h,s}=F_s-F_h
\qquad\text{almost everywhere}\,.
\end{equation}
For every nonnegative measurable function $f\le B_{h,s}$,
\begin{equation}\label{eq:source-energy-bound}
\begin{aligned}
\sum_{k=0}^\infty\langle f,Q_s^kf\rangle
&\le E(s)-E(h),\\
\sum_{k=0}^\infty Q_s^kf
&\le F_s-F_h
\qquad\text{almost everywhere}\,.
\end{aligned}
\end{equation}
In particular, $f\in L^2(\mu)$, and for every integer $T\ge1$,
\begin{equation}\label{eq:finite-sum-L2}
\left\|\sum_{k=0}^{T-1}Q_s^kf\right\|_2^2
\le T[E(s)-E(h)]\,.
\end{equation}
\end{lemma}

\begin{proof}
\emph{Finiteness and local integrability.}
Let $Q_{\mathbb{R}_+^n}$ be defined by~\eqref{eq:Q-definition} with both strip Green operators replaced by $G_{\mathbb{R}_+^n}$. The factorization used in Lemma~\ref{lem:operator}, together with~\eqref{eq:green-contraction}, shows that $Q_{\mathbb{R}_+^n}$ is a contraction on $L^2(\mu)$. Domain monotonicity gives a nonnegative kernel for $Q_{\mathbb{R}_+^n}-Q_h$.

For $x,y\in\supp g_0$, restrict the intermediate variable in the kernel formula to a fixed ball compactly contained in $\{x_n>h\}$. This contribution is absent from $Q_h$. Positivity of the half-space Green kernel and the bounds on the coefficients on compact sets give a positive lower bound for the resulting kernel density with respect to $\mu$. Since $g_0$ is bounded,
\begin{equation*}
(Q_{\mathbb{R}_+^n}-Q_h)(x,\md y)
\ge\kappa_h g_0(x)g_0(y)\md\mu(y)
\end{equation*}
for some $\kappa_h>0$, depending on the fixed solution, $g_0$, and $h$. Outside the product of the supports, the right-hand side is zero.

For an integer $T\ge1$, set $\xi_T=\sum_{k=0}^{T-1}Q_h^kg_0$ and $e_T=\langle g_0,\xi_T\rangle$. Each $\xi_T$ belongs to $L^2(\mu)$ and is nonnegative. Applying the kernel lower bound to $\xi_T$ and using contractivity of $Q_{\mathbb{R}_+^n}$ gives
\begin{align*}
\kappa_h e_T^2
&\le\langle\xi_T,(Q_{\mathbb{R}_+^n}-Q_h)\xi_T\rangle\\
&\le\langle\xi_T,(I-Q_h)\xi_T\rangle\\
&=\langle\xi_T,g_0-Q_h^Tg_0\rangle
\le e_T\,.
\end{align*}
The last inequality uses $Q_h^Tg_0\ge0$. Since $e_T>0$, we obtain $e_T\le\kappa_h^{-1}$, and monotone convergence yields $E(h)\le\kappa_h^{-1}<\infty$.

Kernel symmetry and monotone convergence also give
\begin{equation*}
\langle Q_hg_0,F_h\rangle
=E(h)-\|g_0\|_2^2<\infty\,.
\end{equation*}
On every compact subset of $\Sigma_h$, the function $Q_hg_0$ has a positive lower bound: in its kernel formula, restrict the two integration variables to fixed interior balls, choosing one ball where $g_0>0$. Thus $F_h$ is locally integrable with respect to $\mu$. Since the density $pv^{p+1}$ is bounded above and below by positive constants on compact subsets of $\Sigma_h$, we obtain $F_h\in L^1_{\loc}(\Sigma_h)$ and finiteness almost everywhere.

\emph{Changing the strip height.}
For every integer $m\ge1$,
\begin{equation*}
Q_s^m-Q_h^m
=\sum_{j=0}^{m-1}Q_s^{m-1-j}(Q_s-Q_h)Q_h^j\,.
\end{equation*}
Every term on the right has a nonnegative kernel. Apply this identity to $g_0$, sum over $m$, and use Tonelli's theorem to obtain~\eqref{eq:resolvent-difference}. Both $F_s$ and $F_h$ are finite almost everywhere, so their difference is well defined there.

Kernel symmetry and Tonelli's theorem then give
\begin{align*}
E(s)-E(h)
&=\langle g_0,F_s-F_h\rangle\\
&=\langle F_s,B_{h,s}\rangle\\
&=\langle F_h,B_{h,s}\rangle
+\sum_{k=0}^\infty
\langle B_{h,s},Q_s^kB_{h,s}\rangle\,.
\end{align*}
The left-hand side is finite, and all terms on the right are nonnegative.

\emph{Bounds for the new sources.}
If $0\le f\le B_{h,s}$, positivity of the kernels of $Q_s^k$ gives both inequalities in~\eqref{eq:source-energy-bound}. The $k=0$ term gives $f\in L^2(\mu)$. By self-adjointness of $Q_s$,
\begin{align*}
\left\|\sum_{k=0}^{T-1}Q_s^kf\right\|_2^2
&=\sum_{i,j=0}^{T-1}\langle f,Q_s^{i+j}f\rangle\\
&\le T\sum_{k=0}^\infty\langle f,Q_s^kf\rangle\\
&\le T[E(s)-E(h)]\,,
\end{align*}
because, for each $k$, at most $T$ pairs $(i,j)$ satisfy $i+j=k$. This proves~\eqref{eq:finite-sum-L2}.
\end{proof}

To estimate the total mass of a finite sum, we apply Cauchy--Schwarz inside a horizontal cylinder and estimate the mass outside it. The interior estimate uses the polynomial growth of $\mu$ supplied by Lemma~\ref{lem:energy}.

Choose $\sigma>0$ so that, for every fixed $h>0$ and all $R>0$,
\begin{align}\label{eq:polynomial-measure}
    \mu(B'_R(0)\times(0,h))\le C_h(1+R)^\sigma.
\end{align}
The choice $\sigma=\max\{1,n-\alpha-\beta-2\}$ suffices. To estimate the mass outside the cylinder, we adapt Davies' exponential conjugation argument~\cite[Section 3.2, pp.~83--88]{Davies1989} to finite powers of $Q_h$. The resulting estimates apply uniformly to the source families used below, including after a change of strip height.

\begin{lemma}\label{lem:localization}
Fix $h>0$. For each positive integer $N$, let $f_N\ge0$ be supported in $\Sigma_h\cap\{|x'|\le A\sqrt{N\log(N+2)}\}$ and satisfy $\sum_{k\ge0}\langle f_N,Q_h^kf_N\rangle\le e$, with $A,e$ independent of $N$. Let $H_{N,T}=\sum_{k<T}Q_h^kf_N$. For any fixed $s\ge h$ and $J>0$, there is $B>A$, independent of $N,T$, such that for all sufficiently large integers $N$ and every integer $T$ with $1\le T\le2N$,
\begin{align}\label{eq:localized-tail}
    \int_{\{|x'|>B\sqrt{N\log(N+2)}\}}
    \bigl(H_{N,T}+Q_sH_{N,T}\bigr)\md\mu\le N^{-J}.
\end{align}
For the same $N,T$, the total mass satisfies
\begin{align}\label{eq:mass-upper}
    \int_{\mathbb{R}_+^n} H_{N,T}\md\mu
    \le C N^{1/2+\sigma/4}(\log(N+2))^{\sigma/4}+N^{-J}.
\end{align}
The constants may depend on $n,p,q$, the fixed solution $(u,v)$, the strip heights $h,s$, and $A,e,J$, but are independent of $N$ and $T$.
\end{lemma}

\begin{proof}
We use the kernel and spectral properties in Lemma~\ref{lem:operator} to estimate finite powers with an exponential weight. For $|\varepsilon|h\le a_n/4$, where $a_n$ is given by~\eqref{eq:exponential-row}, define
\[
Q_{h,\varepsilon}(x,\md y)=e^{\varepsilon(|x'|-|y'|)}Q_h(x,\md y).
\]
The inequality $\bigl||x'|-|y'|\bigr|\le|x'-y'|$, the exponential moment bound~\eqref{eq:exponential-row}, and kernel symmetry give boundedness on $L^2(\mu)$ by Schur's test. Decompose $Q_{h,\varepsilon}-Q_h$ into its self-adjoint part $A_\varepsilon$ and skew-adjoint part $B_\varepsilon$. Their kernel factors are $\cosh z-1$ and $\sinh z$, where $z=\varepsilon(|x'|-|y'|)$. Applying the same kernel estimates to these factors gives
\[
A_\varepsilon^*=A_\varepsilon,\quad B_\varepsilon^*=-B_\varepsilon,\qquad
\|A_\varepsilon\|\le C_h\varepsilon^2,\quad \|B_\varepsilon\|\le C_h|\varepsilon|.
\]
The spectral bound $0\le Q_h\le I$ gives $\|Q_hf\|_2^2+\|(I-Q_h)f\|_2^2\le\|f\|_2^2$. When we expand $\|(Q_h+A_\varepsilon+B_\varepsilon)f\|_2^2$, skew-adjointness gives $\operatorname{Re}\langle f,B_\varepsilon f\rangle=0$. The mixed term with $B_\varepsilon$ therefore contains $\operatorname{Re}\langle(Q_h-I)f,B_\varepsilon f\rangle$. Young's inequality bounds this term by half of $\|(I-Q_h)f\|_2^2$ plus $C_h\varepsilon^2\|f\|_2^2$. The terms containing $A_\varepsilon$ and the squares of the perturbations are also bounded by $C_h\varepsilon^2\|f\|_2^2$. Combining these estimates yields
\begin{align}\label{eq:Davies}
    \|Q_{h,\varepsilon}\|_{2\to2}\le e^{C_h\varepsilon^2}.
\end{align}
Consequently, for every integer $k\ge0$, the operator norm of $Q_{h,\varepsilon}^k$ is at most $e^{C_hk\varepsilon^2}$.

Fix $\varepsilon>0$ with $\varepsilon s\le a_n/4$, so that~\eqref{eq:Davies} applies at both heights $h$ and $s$. In a finite product of weighted kernels, the factors at intermediate points cancel, leaving $e^{\varepsilon(|x'|-|y'|)}$. The $k=0$ term in the energy hypothesis gives $\|f_N\|_2\le\sqrt e$. The support condition and~\eqref{eq:Davies} therefore yield
\[
\|e^{\varepsilon|x'|}H_{N,T}\|_2
\le 2N\sqrt e\,\exp\!\left(\varepsilon A\sqrt{N\log(N+2)}+2C_hN\varepsilon^2\right).
\]
Applying $Q_s$ adds only the factor $e^{C_s\varepsilon^2}$. These weighted identities follow from finite products of nonnegative kernels and Tonelli's theorem.

Integrating the radial exponential weight against the volume bound~\eqref{eq:polynomial-measure} gives
\[
\int_{\Sigma_s\cap\{|x'|>R\}}e^{-2\varepsilon|x'|}\md\mu
\le C_s\varepsilon^{-\sigma}e^{-\varepsilon R}\qquad(0<\varepsilon\le1).
\]
Apply Cauchy--Schwarz outside the cylinder of radius $B\sqrt{N\log(N+2)}$, using the two exponential weights above. Take $\varepsilon=\sqrt{\log(N+2)/N}$. Then $N\varepsilon^2$ and $\varepsilon\sqrt{N\log(N+2)}$ both equal $\log(N+2)$. A sufficiently large $B$, independent of $N,T$, therefore gives~\eqref{eq:localized-tail} for both $H_{N,T}$ and $Q_sH_{N,T}$.

The energy hypothesis and the finite-sum calculation in Lemma~\ref{lem:source-energy} give $\|H_{N,T}\|_2^2\le Te\le2Ne$. Inside the cylinder, Cauchy--Schwarz bounds the mass of $H_{N,T}$ by its $L^2$ norm times the square root of the cylinder's measure. Inserting~\eqref{eq:polynomial-measure} and adding the tail estimate proves~\eqref{eq:mass-upper}.
\end{proof}

We now derive lower mass bounds from~\eqref{eq:nonlinear-feedback}. The next lemma also bounds the mean of $\tau=x_n/v$, which will let us choose the next strip height independently of $N$, and selects a partial sum whose first omitted iterate has small mass relative to the sum. This last estimate will control the residual when the height is increased.

Fix $h>0$ and a nonnegative source $f\in L^1(\mu)$ supported in $\Sigma_h$. For integers $k\ge0$ and $T\ge1$, set
\[
g_k=Q_h^kf,\qquad m_k=\int_{\mathbb{R}_+^n}g_k\md\mu,\qquad
q_k=\int_{\mathbb{R}_+^n}r_hg_k\md\mu,\qquad
M_T=\sum_{k<T}m_k.
\]
Thus $M_T$ is the total mass of the partial sum, and $q_k$ is the mass of $g_k$ weighted by $1-N_h/v$ on $\Sigma_h$.

\begin{lemma}\label{lem:finite-sums}
All $m_k,q_k,M_T$ are finite. If $q_0>0$, then $q_k>0$ for all $k\ge0$ and, for every integer $T\ge1$,
\begin{align}\label{eq:quadratic-mass}
    M_T\ge\delta_*q_0(T-1)^2.
\end{align}
If $q_0>0$ and $a_0:=\int_{\mathbb R_+^n}\tau f\md\mu<\infty$, define $A_T:=\int_{\mathbb R_+^n}\tau\sum_{k<T}g_k\md\mu$. Then $A_T<\infty$ for every integer $T\ge1$, and
\begin{align}\label{eq:moment-recursion}
    \frac{A_T}{M_T}\le\frac{a_0}{M_T}
    +K_h\left(\frac{A_T}{M_T}\right)^{1/(pq)}.
\end{align}
For a family $f_N$ satisfying the support and uniform energy hypotheses of Lemma~\ref{lem:localization}, define $m_{k,N},M_{T,N}$ as above. If $M_{N,N}\ge c_*>0$ uniformly in $N$, then some $T_N\in[N,2N-1]\cap\mathbb Z$ satisfies
\begin{align}\label{eq:terminal-ratio}
    \frac{m_{T_N,N}}{M_{T_N,N}}\le C\frac{\log(N+2)}N\,,
\end{align}
where $C$ is independent of $N$.
\end{lemma}

\begin{proof}
Kernel symmetry gives $m_{k+1}=\langle Q_h1,g_k\rangle\le pq\,m_k$, so $m_k,q_k,M_T$ are finite. Assume $q_0>0$. Pairing~\eqref{eq:nonlinear-feedback} with $g_k$ gives the first inequality below and shows inductively that $q_{k+1}\ge q_k>0$. Hence $m_{k+1}>0$, and Cauchy--Schwarz gives the second inequality:
\[
q_{k+1}-q_k\ge\delta_*\int_{\mathbb{R}_+^n}r_h^2g_{k+1}\md\mu
\ge\delta_*\frac{q_{k+1}^2}{m_{k+1}}.
\]
Dividing by $q_kq_{k+1}$ and summing over $k$ yield
\begin{align*}
    \begin{aligned}
        \frac1{q_k}-\frac1{q_{k+1}}&\ge\frac{\delta_*}{m_{k+1}},\\
        \sum_{k=1}^{T-1}\frac1{m_k}&\le\frac1{\delta_*q_0}.
    \end{aligned}
\end{align*}
Cauchy--Schwarz applied to $\sum_{k=1}^{T-1}1$ now yields~\eqref{eq:quadratic-mass}. The source enters this bound only through $q_0$, so the estimate can be applied separately to each member of a varying source family.

If $a_0<\infty$,~\eqref{eq:Foster} and H\"older's inequality show that a finite mass and a finite $\tau$-moment of $g_k$ give a finite $\tau$-moment of $g_{k+1}$. Summing these moment estimates and applying H\"older's inequality gives~\eqref{eq:moment-recursion}. Since $1/(pq)<1$, a uniform bound on $a_0/M_T$ then bounds $A_T/M_T$ uniformly.

For the last conclusion, fix $N$ and use the same source $f_N$ for all partial sums with indices between $N$ and $2N$. The upper bound~\eqref{eq:mass-upper} and the assumption $M_{N,N}\ge c_*$ give
\[
\sum_{T=N}^{2N-1}\log\left(1+\frac{m_{T,N}}{M_{T,N}}\right)
=\log\frac{M_{2N,N}}{M_{N,N}}\le C\log(N+2)\,.
\]
The $N$ logarithmic increments are nonnegative, so at least one is at most $C\log(N+2)/N$. For sufficiently large $N$, this gives~\eqref{eq:terminal-ratio}.
\end{proof}

We start the height iteration with the same $g_0$ fixed before Lemma~\ref{lem:source-energy}. The function $\tau$ is bounded on $\supp g_0$, and~\eqref{eq:strip-decay} gives $l_h=N_h/v\to0$ uniformly there as $h\to\infty$. Choose a fixed $h_0$ above the support such that $\langle r_{h_0},g_0\rangle>0$. The quadratic series of $g_0$ is $E(h_0)<\infty$, so its compact support allows us to apply Lemmas~\ref{lem:localization} and~\ref{lem:finite-sums} with $f_N=g_0$. For each sufficiently large $N$, choose $T_{0,N}$ by~\eqref{eq:terminal-ratio} and set
\[
H_{0,N}=\sum_{k<T_{0,N}}Q_{h_0}^kg_0,\qquad
L_{0,N}=\int_{\mathbb{R}_+^n}H_{0,N}\md\mu.
\]
Equations~\eqref{eq:quadratic-mass} and~\eqref{eq:moment-recursion} give $L_{0,N}\ge c_0N^2$ and $L_{0,N}^{-1}\int_{\mathbb{R}_+^n}\tau H_{0,N}\md\mu\le C_0$, with constants independent of $N$. The construction also gives $H_{0,N}\le F_{h_0}$, $\|H_{0,N}\|_2^2\le2NE(h_0)$, and the tail estimate~\eqref{eq:localized-tail}. Our choice of $T_{0,N}$ ensures $\int Q_{h_0}^{T_{0,N}}g_0\md\mu=o(L_{0,N})$. Thus the first omitted iterate has small mass relative to $H_{0,N}$ in the identity $(I-Q_{h_0})H_{0,N}=g_0-Q_{h_0}^{T_{0,N}}g_0$.

These estimates start the induction. The next proposition passes to one fixed larger height, gains a factor of order $N^2$ in mass, and preserves the bounds needed for another height increase.

\begin{proposition}\label{prop:height-transfer}
Fix the source $g_0$ in Lemma~\ref{lem:source-energy} and $h$ above its support. For all sufficiently large integers $N$, let $f_N\ge0$ belong to $L^1(\mu)$, be supported in $\Sigma_h\cap\{|x'|\le A\sqrt{N\log(N+2)}\}$, and satisfy
\[
 \sum_{k\ge0}\langle f_N,Q_h^kf_N\rangle\le e_h.
\]
For $T_N\in[N,2N-1]\cap\mathbb Z$, set $H_N=\sum_{k<T_N}Q_h^kf_N$ and $L_N=\int H_N\md\mu$. Assume
\[
 H_N\le F_h,\qquad L_N\ge c_hN^2,\qquad
 \int\tau H_N\md\mu\le C_{\tau,h}L_N,\qquad
 \int Q_h^{T_N}f_N\md\mu=o(L_N),
\]
where $A,e_h,c_h,C_{\tau,h}$ are independent of $N$. There is a fixed $s>h$ and compactly supported $f'_N\ge0$ in $\Sigma_s$, with horizontal support radius at most $A_s\sqrt{N\log(N+2)}$, such that
\[
 \sum_{k\ge0}\langle f'_N,Q_s^kf'_N\rangle\le E(s)-E(h).
\]
There are $T'_N\in[N,2N-1]\cap\mathbb Z$ for which $H'_N=\sum_{k<T'_N}Q_s^kf'_N$, $L'_N=\int H'_N\md\mu$ satisfy
\[
 H'_N\le F_s,\qquad L'_N\ge c'_hN^2L_N,\qquad
 \int\tau H'_N\md\mu\le C'_{\tau,h}L'_N,\qquad
 \int Q_s^{T'_N}f'_N\md\mu=o(L'_N).
\]
The new height and constants may depend on the fixed solution, $g_0$, the fixed heights and previous bounds, but not on $N$.
\end{proposition}

\begin{proof}
The support and energy assumptions allow us to apply Lemma~\ref{lem:localization} to $H_N+Q_sH_N$ at each fixed $s\ge h$. We first choose the new height using the uniform $\tau$-moment bound. Let $c_{\mathrm{strip}}$ be the constant $c$ in Lemma~\ref{lem:strip-decay}. Choose $a_*\in(0,1)$ such that
\begin{equation*}
pq(a_*^\gamma+a_*^{1/(pq)})\le\frac{pq-1}{4}\,.
\end{equation*}
On $\Sigma_h$,~\eqref{eq:strip-decay} gives $l_s\le Cs^{-\beta-1}\tau$ whenever $h\le c_{\mathrm{strip}}s$. Hence one fixed $s>h$ can be chosen so large that $h\le c_{\mathrm{strip}}s$ and
\[
\frac1{L_N}\int_{\mathbb{R}_+^n}l_sH_N\md\mu
\le Cs^{-\beta-1}\frac1{L_N}\int_{\mathbb{R}_+^n}\tau H_N\md\mu\le a_*.
\]
The choice of $s$ is independent of $N$ because the normalized $\tau$-moments are uniformly bounded.

Apply~\eqref{eq:helper} and Jensen's inequality to the probability measure $H_N\,\md\mu/L_N$. The exponents $\gamma$ and $1/(pq)$ are at most one, so the preceding average bound gives
\[
\langle H_N,(Q_s-I)r_s\rangle
\ge\bigl[pq-1-pq(a_*^\gamma+a_*^{1/(pq)})\bigr]L_N
\ge\frac{3(pq-1)}4L_N.
\]
Define the new source before truncation by $b_N=(Q_s-Q_h)H_N\ge0$. Kernel symmetry and the identity $(I-Q_h)H_N=f_N-Q_h^{T_N}f_N$ give
\begin{align*}
    \begin{aligned}
        \langle r_s,b_N\rangle
        &=\langle H_N,(Q_s-I)r_s\rangle
        +\langle r_s,f_N\rangle-\langle r_s,Q_h^{T_N}f_N\rangle\\
        &\ge\frac{pq-1}{2}L_N
    \end{aligned}
\end{align*}
for sufficiently large $N$. Indeed, $\langle r_s,f_N\rangle\ge0$, while $\langle r_s,Q_h^{T_N}f_N\rangle\le\int Q_h^{T_N}f_N\md\mu=o(L_N)$. All pairings are finite because $0\le r_s\le1$, the kernel masses are bounded by $pq$, and $H_N,f_N\in L^1(\mu)$.

Since $0\le b_N\le Q_sH_N$,~\eqref{eq:localized-tail} at this fixed $s$, with $J=10$, bounds the mass of $b_N$ outside horizontal radius $A_s\sqrt{N\log(N+2)}$ by $N^{-10}$. Since $L_N\ge c_hN^2$, this loss is $o(L_N)$. After this restriction, choose a compact interval of heights inside $(0,s)$ for each $N$. Monotone convergence and the lower bound just proved give a nonnegative compactly supported source $f'_N\le b_N$ such that
\begin{align}\label{eq:retained-work}
    \langle r_s,f'_N\rangle\ge\frac{pq-1}{4}L_N.
\end{align}
Both restrictions act only on $b_N$; the interval of heights is chosen separately for each $N$, after $s$ has been fixed.

The inequality $H_N\le F_h$ and positivity of the kernel of $Q_s-Q_h$ give $f'_N\le b_N\le(Q_s-Q_h)F_h=B_{h,s}$. Thus~\eqref{eq:source-energy-bound} yields
\[
\sum_{k\ge0}\langle f'_N,Q_s^kf'_N\rangle\le E(s)-E(h),\qquad
\sum_{k\ge0}Q_s^kf'_N\le F_s-F_h.
\]
Since $h,s$ and $g_0$ are fixed, this quadratic-series bound is independent of $N$. For every integer $T\in[N,2N]$,~\eqref{eq:quadratic-mass} and~\eqref{eq:retained-work} give
\begin{equation*}
M'_T:=\int_{\mathbb{R}_+^n}\sum_{k<T}Q_s^kf'_N\md\mu
\ge c'N^2L_N\,.
\end{equation*}
The new source also satisfies
\[
\int_{\mathbb{R}_+^n}f'_N\md\mu\le pq\,L_N,\qquad
\int_{\mathbb{R}_+^n}\tau f'_N\md\mu
\le\langle Q_s\tau,H_N\rangle\le K_s C_{\tau,h}^{1/(pq)}L_N.
\]
It follows that $(M'_T)^{-1}\int\tau f'_N\md\mu=O(N^{-2})$ for $N\le T\le2N$. Since $1/(pq)<1$,~\eqref{eq:moment-recursion} then bounds the normalized $\tau$-moments of the new partial sums uniformly in $N,T$.

The new sources satisfy the support and energy hypotheses of Lemma~\ref{lem:localization}. Apply that lemma and choose $T'_N\in[N,2N-1]$ by~\eqref{eq:terminal-ratio}. Set $H'_N=\sum_{k<T'_N}Q_s^kf'_N$. The preceding estimates give $L'_N\ge c'N^2L_N$ and a uniform bound on its normalized $\tau$-moment. Equations~\eqref{eq:source-energy-bound} and~\eqref{eq:finite-sum-L2} give $H'_N\le F_s-F_h\le F_s$ and $\|H'_N\|_2^2\le2N[E(s)-E(h)]$. The choice of $T'_N$ gives $\int Q_s^{T'_N}f'_N\md\mu=o(L'_N)$. Thus all conditions needed for another height increase hold at $s$, with the same fixed $g_0$ supplying $F_s$ and $E(s)$.
\end{proof}

We now complete the proof using Proposition~\ref{prop:height-transfer} and~\eqref{eq:mass-upper}.

\begin{proof}[Proof of Theorem~\ref{thm:main}]
A nonnegative, nonzero solution has two positive components, as noted in the introduction. If $n=1$, both functions are strictly concave on $(0,\infty)$. Positivity on the entire half-line forces their first derivatives to be nonnegative. Then $v(x_n)\ge v(1)>0$ for $x_n\ge1$, so $u''\le-v(1)^p$ forces $u'$ eventually to be negative, a contradiction.

Let $n\ge2$. Lemma~\ref{lem:small-product} handles $pq\le1$. In the remaining case, interchange the components and exponents if necessary so that $p\ge q>0$, and use Lemma~\ref{lem:global-reduction} to obtain a positive solution stable on all of $\mathbb{R}_+^n$. Apply the construction of this section to that one solution.

Fix an integer $K\ge0$ such that $2+2K>1/2+\sigma/4$. Repeating Proposition~\ref{prop:height-transfer} a finite number of times gives fixed heights $h_0<\cdots<h_K$ and, for every sufficiently large $N$, a final sum with
\[
L_{K,N}\ge c_KN^{2+2K},\qquad
\|H_{K,N}\|_2\le C_K\sqrt N.
\]
Its horizontal tail is bounded by $N^{-10}$ outside a cylinder of radius $C_K\sqrt{N\log(N+2)}$. The growth estimate~\eqref{eq:polynomial-measure} and Cauchy--Schwarz therefore imply
\[
L_{K,N}\le C_KN^{1/2+\sigma/4}(\log(N+2))^{\sigma/4}+N^{-10}.
\]
The strict inequality defining $K$ makes these two bounds incompatible as $N\to\infty$.
\end{proof}

\section*{Declarations}

\noindent\textbf{Conflict of interest.}
The authors declare that they have no conflicts of interest.

\medskip

\noindent\textbf{Data availability.}
Data sharing is not applicable to this article, as no datasets were generated or analyzed during the current study.

\medskip

\noindent\textbf{Generative AI disclosure.}
ChatGPT was used to assist in discussing proof strategies, verifying mathematical arguments, reviewing relevant references, and preparing the manuscript. The authors have reviewed all AI-assisted material, further refined the exposition, and take full responsibility for the mathematical results and the final version of the manuscript.

\bibliographystyle{abbrvurl}
\bibliography{references}

\end{document}